\documentclass[a4paper, 11pt]{article}

\usepackage{longtable}
\usepackage{geometry}
\usepackage{graphicx}
\usepackage{subfig}
\usepackage{color}
\usepackage{xcolor}
\usepackage{bm}
\usepackage{fixmath}
\usepackage{mathrsfs}
\usepackage{ulem}
\usepackage{multirow}
\usepackage{pdflscape}
\usepackage{rotating}
\usepackage{lscape}
\usepackage{pdflscape}
\usepackage{titlesec}
\usepackage{float}
\usepackage{hyperref}\usepackage[utf8]{inputenc}
\usepackage[english]{babel}
\usepackage{amsthm}
\usepackage{multirow}

\begin{document}







\def\b{\begin{eqnarray}}
\def\e{\end{eqnarray}}
\def\n{\noindent}

%


\begin{center}
{\huge \textbf{Isolation Intervals of the Real Roots of \vskip.35cm  the Parametric Cubic Equation and  \vskip.5cm Improved Complete Root Classification}}


\vspace{9mm}
\noindent
{\Large \bf Emil M. Prodanov} \vskip.7cm
{\it School of Mathematical Sciences, Technological University Dublin,
\vskip.1cm
Park House, Grangegorman, 191 North Circular Road, Dublin D07 EWV4, Ireland,}
\vskip.1cm
{\it E-Mail: emil.prodanov@tudublin.ie} \\
\vskip.5cm
\end{center}

\vskip1cm


\begin{abstract}
\n
The isolation intervals of the real roots of the real symbolic monic cubic polynomial $p(x) = x^3 + a x^2 + b x + c$ are found in terms of  simple functions of the coefficients of the polynomial (such as: $-a$, $-a/3$, $-c/b$, $\pm \sqrt{-b}$, when $b$ is negative), and the roots of some {\it auxiliary quadratic equations} whose coefficients are also simple functions of the coefficients of the cubic. All possible cases are presented with clear and very detailed diagrams. It is very easy to identify which of these diagrams is the relevant one for any given cubic equation and to read from it the isolation intervals of the real roots of the equation. A much-improved complete root classification, addressing the signs (together with giving the isolation intervals) of the individual roots, is also presented. No numerical approximations or root finding techniques are used. Instead of considering the discriminant of the cubic, criterion for the existence of a single real root or three real roots is found as conditions on the coefficients of the cubic, resulting from the roots of the auxiliary quadratic equations. It is also shown that, if a cubic equation has three real roots, then these lie in an interval $I$ such that $\sqrt{3}\sqrt{a^2/3 - b} \le I \le 2 \sqrt{a^2/3 - b}$, independent of $c$. A detailed algorithm for applying the method for isolation of the roots of the cubic is also given and it is illustrated through examples, including the full mathematical analysis of the cubic equation associated with the Rayleigh elastic waves and finding the isolation intervals of its real roots.
\end{abstract}



\newgeometry{top=35mm, bottom=35mm,left=35mm, right=35mm}
\section{Introduction}

The explicit formul\ae \, for the roots of the cubic equation were discovered during the Renaissance \cite{card} and since then, there has been significant research on the theory of equations --- see, for example, the reading list \cite{g} from 1933. For more recent and comprehensive reference, see \cite{rs}--\cite{k}. \\
The general cubic equation $Ax^3 + Bx^2 + Cx + D = 0$ can be {\it depressed} by the coordinate translation $x \to x - B/3A$ and dividing by $A$ afterwards to obtain the depressed cubic equation $x^3 + px + q = 0$, for which $p = (3AC - B^2)/(3A^2)$ and $q = (2 B^3 - 9 ABC + 27 A^2 D)/(27 A^3)$. The Cardano cubic formul\ae \, apply to the {\it depressed cubic equation} and yield its roots. The roots of the original equation can then be recovered by making the inverse coordinate transformation.  \\
Despite of the existence of the explicit Cardano cubic formul\ae, it is not always easy to apply them. For example, when the roots of a real cubic equation are all real and distinct, the Cardano formul\ae \, give the roots in a form involving cube roots of complex numbers and the cube root of a general complex number cannot be expressed in the form $a + bi$, where $a$ and $b$ involve only real radicals \cite{d}. This is the so called {\it casus irreducibilis}. The cube root of a complex number appears in result of the necessity to introduce an imaginary number in the Cardano formul\ae \, by taking the square root of a negative number.  Secondly, when the coefficients of the equation are parameters or functions of some parameters, the above renders the Cardano formul\ae \, practically inapplicable. \\
The first remedy in such situations would be to study the root locations, that is, to find intervals which contain the roots of the equation. If a finite interval is found such that only one root of the equation lies in it, then this interval is called {\it isolation interval}.  \\
The Descartes Rule of Signs \cite{c} from 1637 was the first step in the direction of locating the roots of an equation. This rule yields the isolation of more than one root over $x > 0$ (or over $x < 0$, when $x$ is replaced by $-x$ in the equation) with some indeterminacy: {\it the number of positive real roots of an equation with real coefficients cannot be greater than the number of the variations of sign in the sequence of its coefficients} (Descartes' original formulation). In 1876, Gauss \cite{ga} found a more precise formulation: {\it the number of positive real roots (counted with their multiplicities) is equal to the number of the variations of sign or is less than that number by a positive even integer}. Newton’s rule \cite{newton}, giving an upper limit on the number of positive and negative roots, could also be hard to apply due to need to analyse the signs of the associated parametric quadratic elements. \\
There has also been an extensive amount of research in another direction: on the determination of the number of roots of an equation over some prescribed interval. For example, this can be done using Sturm's \cite{st} or Budan's \cite{bu} theorems, see also \cite{d}. \\
The Sturm sequence for the monic cubic polynomial $p(x) = x^3 + a x^2 + b x + c$ is
\b
p_0(x) & = & p(x) \, \, = \,\, x^3 + a x^2 + b x + c, \\
p_1(x) & = & p'(x) \,\, = \,\, 3 x^2 + 2 a x + b, \\
p_2(x) & = & - \mbox{rem} \left[ \frac{p_0(x)}{p_1(x)} \right] \,\, = \,\, \frac{2}{3} \left( \frac{a^2}{3} - b \right) x + \frac{2}{9} ab - c, \\
p_3 & = & - \mbox{rem} \left[ \frac{p_1(x)}{p_2(x)} \right] \,\, = \,\, - \frac{1}{4} \frac{\Delta_{(3)}}{(\frac{a^2}{3} - b)^2} \,\, = \,\, \mbox{const},
\e
where rem$[p_i(x)/p_{i+1}(x)]$ denotes the remainder of the division of the polynomial $p_i(x)$ by the polynomial $p_{i+1}(x)$ and
\b
\label{dis}
\Delta_{(3)} = -27 c^2 + (18 a b - 4 a^3) c + a^2 b^2 - 4b^3
\e
is the discriminant of $p(x)$. [The product of the squares of the differences of the roots $x_i$ of a monic polynomial is called {\it discriminant} of the polynomial. For the cubic, it is $\Delta_{(3)} = (x_1 - x_2)^2 \, (x_2 - x_3)^2 \, (x_3 - x_1)^2$.] \\
Sturm's theorem \cite{d} allows the determination of the number of real roots of $p(x) = 0$ between $a$ and $b$. This number is equal to the excess of the number of variations of sign of $p(x), \, p_1(x), \, p_2(x),$ and $p_3$ for $x = a$ over the number of variations of sign of these for $x = b$ (the vanishing terms are ignored). \\
Cheng and Lin \cite{cheng} consider a nonlinear function with $n$ parameters and determine if none, some or all of the roots of the function lie in a specified subregion of the domain of the function for a partition of the parameter region. \\
Another vein in the research on the theory of equations are the {\it root classification} and the {\it complete root classification} of parametric polynomials. These refer to the provision of all possible cases of the polynomial roots and consist of the list of the multiplicities of all roots, in the case of {\it root classification}, and the root classification, together with the conditions which the equation coefficients should satisfy for each case of the root classification, in the case of {\it complete root classification}. The {\it root classification} and the {\it complete root classification} do not determine the location of the roots. \\
For the very simple case of a depressed cubic polynomial $x^3 + px + q$, the {\it complete root classification} is \cite{d}, \cite{a}:
\begin{center}
\begin{tabular}{|l|r|}
\hline
$\delta_{(3)} > 0$ &  $\{ 1, 1, 1 \}$ \\
\hline
$\delta_{(3)} = 0$ &  $\{ 1, 2 \}$  \\
\hline
$\delta_{(3)} < 0$ &  $\{ 1 \}$ \\
\hline
\end{tabular}
\end{center}
where $\delta_{(3)} = -4 p^3 - 27 q^2$ is the discriminant of $x^3 + px + q$ and the lists in the figure brackets give the real root multiplicities. \\
This paper proposes a method with which the isolation intervals of the roots of the symbolic monic cubic polynomial $p(x) = x^3 + a x^2 + b x + c$ are found. The criteria for existence of a single real root or three real roots are found as conditions on the coefficients of the cubic equation resulting from the roots of some {\it auxiliary quadratic equations}. The endpoints of these isolation intervals are found as roots of the {\it auxiliary quadratic equations} or simple functions of the coefficients of the cubic: $-a$, $-a/3$, $-c/b$, $\pm \sqrt{-b}$ (when $b$ is negative), etc. \\
Also presented in this paper is a much-improved {\it complete root classification}, giving the conditions (as well as the isolation intervals) for three positive roots, for two positive and one negative root, for one positive and two negative roots, for three negative roots, for a single positive root, and for a single negative root. \\
All results are precise --- no numerical approximations or root finding algorithms have been used. \\
One of the examples given in this work provides a thorough analysis of the parametric cubic equation associated with the Rayleigh elastic waves. The isolation intervals of its real roots are found for different values of the parameter of the equation. This illustrates how the proposed method allows one to analyze parametric cubic equations which are a very common occurrence in all branches of mathematics, sciences and engineering. \\
The isolation intervals of the real roots of the general parametric quartic equation and its complete root classification have been studied in \cite{3} and those of the quintic equation --- in \cite{28}, and, with diminishing degree of determinacy, recursively for the general polynomial, --- in \cite{1}.


\section{The Auxiliary Quadratic Equations}
\n
Consider the general monic cubic polynomial
\b
\label{p}
p(x) = x^3 + a x^2 + b x + c.
\e
Its discriminant $\Delta_{(3)} = -27 c^2 + (18 a b - 4 a^3) c + a^2 b^2 - 4b^3$ is quadratic in the free term $c$ of the polynomial (as $\delta_{(3)}$ is in $q$). In turn, the discriminant of this quadratic is
\b
\Delta_{(2)} = 16 (a^2 - 3b)^3.
\e
Given that the leading coefficient of $\Delta_{(3)}$ is negative, if $b > a^2/3$, then $\Delta_{(2)} < 0$ for all $a$. Therefore, $\Delta_{(3)} < 0$ for all $a$, all $b > a^2/3$, and all $c$. Hence, the polynomial $x^3 + a x^2 + b x + c$ with $b > a^2/3$ for any $a$ and for any $c$ will have only one real root (and two complex conjugate roots). \\
The monic cubic polynomial $x^3 + a x^2 + b x + c$ with $b < a^2/3$ may have either one real root or three real roots (they are distinct if $\Delta_{(3)} > 0$ and there is a double root if $\Delta_{(3)} = 0$). To determine which of these occurs, one needs to solve the inequality $\Delta_{(3)} \ge 0$. \\
Firstly, consider the quadratic equation $\Delta_{(3)} = 0$, that is
\b
\label{q1}
c^2 + \left( \frac{4}{27} \, a^3 - \frac{2}{3}\, a b \right) c  - \frac{1}{27} \, a^2 b^2  + \frac{4}{27}\, b^3 = 0.
\e
This is the {\it first auxiliary quadratic equation}. Its roots are:
\b
\label{c12}
c_{1,2}(a,b) = c_0 \, \pm \, \frac{2}{27} \, \sqrt{( a^2 - 3 b)^3},
\e
where
\b
\label{c0}
c_0(a,b) = -\frac{2}{27} \, a^3 + \frac{1}{3} \, a b.
\e
Clearly, for any $a$ and $b < a^2/3$, when $c_2 < c < c_1$, the cubic polynomial $x^3 + a x^2 + b x + c$ will have three distinct real roots (as its discriminant $\Delta_{(3)}$ will be positive). \\
For any $a$ and $b  < a^2/3$, when $c = c_1$ or $c = c_2$, the discriminant $\Delta_{(3)}$ of the cubic polynomial $x^3 + a x^2 + b x + c$ will be zero and, hence, the polynomial will have three real roots, two of which equal. The roots in this case can be easily determined using the Vi\`ete formul\ae \,\ for the cubic equation $x^3 + a x^2 + b x + c = 0$ in the case of a double root $\mu_i$ and a simple root $\xi_i$ (with $i = 1$ when $c = c_1$ and $i = 2$ when $c = c_2$), that is, using $2 \mu_i + \xi_i = - a, \,\, \mu_i^2 + 2 \mu_i \xi_i = b,$ and $\mu_i^2 \xi_i = -c_i$ (for $i = 1, 2$).
The first formula gives $\xi_i = - a - 2 \mu_i$. Substituting into the second, leads to
\b
\label{q2}
3 \mu_i^2 + 2 a \mu_i + b = 0.
\e
This is the {\it second auxiliary quadratic equation}. \\
In the regime $b \le a^2/3$, the real roots of equation (\ref{q2}), namely, the double root $\mu_{1,2}$ of the cubic with $c = c_{1,2}$, are
\b
\label{mu12}
\mu_{1,2} = -\frac{a}{3} \pm \frac{\sqrt{3}}{3} \sqrt{\frac{a^2}{3} - b}.
\e
The corresponding simple real root of the cubic with $c = c_{1,2}$ is
\b
\label{xi}
\xi_{1,2}  =  - a - 2 \mu_{1,2} = -\frac{a}{3} \mp \frac{2\sqrt{3}}{3} \sqrt{\frac{a^2}{3} - b}.
\e
Finally, for any $a$ and $b = a^2/3$, one has $c_1 = c_2 = a^3/27$. The cubic polynomial with $c = a^3/27$ will be exactly $(x + a/3)^3$ and it will have a triple real root $-a/3$. For any other value of $c$, when $b = a^2/3$, one has $\Delta_{(3)} = - (a^3 - 27 c)^2/27$, which is negative. The cubic in this case will have only one real root $- a/3 + \sqrt[3]{a^3/27 - c}\,$ (found by completing the cube). \\
Note that equation (\ref{q2}) is nothing else but the equation for the critical points of the cubic polynomial $x^3 + a x^2 + b x + c$. The two ``extreme" cubics, namely $x^3 + a x^2 + b x + c_1$ and $x^3 + a x^2 + b x + c_2$, are such that the graph of each of them is tangent to the abscissa at the double root, that is, for the ``extreme" cubics the local extrema (the maximum of the cubic with $c = c_2$ and the minimum of the cubic $c = c_1$, both for $b < a^2/3$) and, also, the saddle of the cubic with $b = a^2/3$ and $c = a^3/27$ lie on the abscissa. \\
Note also that the free terms of the two ``extreme" cubics are $c_i = c - p(\mu_i)$. \\
The cubic with $c = c_0 = - 2 a^3/27 + a b/3$ has three real roots. These are $\rho_0 = - a/3$ (which is the first coordinate projection of the inflection point of this particular cubic as well as of the general one) and $\rho_{1,2}$ which are equidistant from $\rho_0$ and are given by
\b
\label{rho}
\rho_{1,2} = -\frac{a}{3} \pm \sqrt{\frac{a^2}{3} - b}.
\e
Note that the critical points $\mu_{1,2} = - a/3 \pm (\sqrt{3}/3) \sqrt{a^2/3-b}$ of the cubic polynomial are equidistant from its inflection point. \\
In the case of a cubic equation with three real roots, the length of the interval where all roots are, varies between its minimum value, achieved when $c = c_i \, (i=1,2)$, namely between $ | \mu_i - \xi_i | = \sqrt{3} \sqrt{a^2/3 - b}$ and its maximum value, $2\sqrt{a^2/3 - b}$ achieved when $c = c_0$, see also \cite{4}. That is, if a cubic equation has three real roots, then the length of the interval $I$ which contains all three roots satisfies
\b
\label{harness}
\sqrt{3} \, \sqrt{\frac{a^2}{3} - b} \le I \le 2 \sqrt{\frac{a^2}{3} - b}
\e
--- independent of the equation parameter $c$. This ``root harness" is another constraint on the roots of the cubic and works in conjunction with the isolation intervals of the roots.

\section{The Method}
\n
The essence of the method, based on the ideas of \cite{1}, is to re-write the cubic equation $x^3 + a x^2 + b x + c = 0$ as
\b
\label{split}
x^2 (x + a) = - b x - c
\e
and seek the intersection points of the ``sub-cubic" $x^2 (x + a)$ with the straight line $- b x - c$. An important feature of $x^2 (x + a)$ is that it has a simple root at $-a$ and a double root at zero. The analysis of the intersection points, that is, the real roots of the cubic equation, is done by studying $x^2 (x + a)$ in specific ranges of $a$ and the straight line $- b x - c$ in specific ranges of $b$, while allowing the variation of the free term $c$, namely, allowing the straight line to ``slide" vertically and, in this process, reveal the different root scenarios. \\
The isolation intervals of the real roots of the cubic can be read graphically by analysing the position of the roots relative to a number of fixed points from a set of five suitable parallel straight lines. These straight lines are determined as follows. Firstly, all these straight lines have the same slope $-b$ as the straight line $ - b x - c$ on the right-hand side of (\ref{split}). \\
Two of these straight lines are the ones which are tangent to the graph of the left-hand side $x^2(x+a)$ of (\ref{split}). They are given by $- b x - c_{1,2}$ and they intersect the ordinate at $-c_{1,2}$. \\
Another straight line of this set is the one that goes through the inflection point of the cubic which coincides with the inflection point of the left-hand side $x^2(x+a)$. This is the straight line $- b x - c_0$, intersecting the ordinate at $-c_0$. Note that $- c_1 \le - c_0 \le - c_2$. \\ Another straight line is $- b x - ab$. This line passes through the point at which $x^2 (x + a)$ crosses the abscissa and through point $-ab$ from the ordinate. Depending on the values of $a$ and $b$, in view of $c_0  = - 2 a^3/27 + a b/3$, one can have $-c_0 < - ab$ or $c_0 \ge -ab$ and one can also have $-c_0 < 0$ or $-c_0 \ge 0$. \\
The final straight line is the separatrix $-bx$. This corresponds to a cubic with $c = 0$, that is $x^3 + a x^2 + b x$. Such cubic has a zero root and two more roots which are given by the roots
\b
\label{lala}
\lambda_{1,2} = -\frac{a}{2} \pm \sqrt{\frac{a^2}{4} - b}
\e
of the {\it third auxiliary quadratic equation}:
\b
\label{q3}
x^2 + a x + b = 0.
\e
Note that if $b > a^2/4$, the roots $\lambda_{1,2}$ are not real. This means that $x^2(x + a)$ intersects the straight line $-bx$ only once --- at the origin. Alternatively, if $b \le a^2/4$, then $\lambda_{1,2}$ are both real and hence $x^2(x + a)$ intersects the straight line $-bx$, except at the origin, at two more points: one in the first quadrant and one in the third, should $b < 0$, and one in the second quadrant and one in the fourth, should $b > 0$ (if $b = 0$, the straight line $- bx$ is the abscissa itself and it intersects $x^2(x + a)$ at $-a$ and $0$). Hence, if $\lambda_{1,2}$ are not real, that is, if $b > a^2/4$, then $c_1$ and $c_2$ will have the same sign. This situation is shown on Figures 14 and 15. Note that $\mu_{1,2}$ ar real only when $b \le a^2/3$ and the straight lines $- b x - c_{1,2}$ exist only in these cases. \\
These are the different regimes for the coefficient $b$:
\begin{itemize}
\item[{\bf (1)}] $\bm{b < - a^2/9}$. Thus, $- ab < - c_0 < 0$ for $a < 0$ (Figure 4) and $0 < -c_0 < - ab$ for $a > 0$ (Figure 5). For both, $\mu_{1,2}$ are real. The cubic can have either three real roots (for $c_2 \le c \le c_1$) or only one real root.
\item[{\bf (2)}] $\bm{- a^2/9 \le b < 0}$. Thus $- c_0 \le - ab < 0$ for $a < 0$ (Figure 6) and $0 < - ab \le - c_0$ for $a > 0$ (Figure 7).  For both, $\mu_{1,2}$ are real. The cubic can have either three real roots (for $c_2 \le c \le c_1$) or only one real root.
\item[{\bf (3)}] $\bm{b = 0}$. See Figure 8 with $a < 0$ and Figure 9 with $a > 0$.  For both, $\mu_{1,2}$ are real ($\mu_2 = 0$ and $\mu_1 = -2a/3$). The cubic can have either three real roots (for $c_2 \le c \le c_1$) or only one real root.
\item[{\bf (4)}] $\bm{0 < b \le 2 a^2/9}$. Thus, if $a < 0$, one has $- c_0 < 0$ and $- ab > 0$ (Figure 10). Alternatively,  if $a > 0$, one has $- c_0 >0$ and $- ab < 0$ (Figure 11). For both, $c_0 = 0$ if $b = 2 a^2/9$. Also for both, $\mu_{1,2}$ are real. The cubic can have either three real roots (for $c_2 \le c \le c_1$) or only one real root.
\item[{\bf (5)}] $\bm{2 a^2/9 < b \le a^2/4}$. Thus $\lambda_{1,2}$ are both real and $c_1$ and $c_2$ have opposite sign: $- c_1 < 0$ and $- c_2 > 0$. Also, $-c_0$ and $-ab$ have the same sign: they are both positive if $a$ is negative (Figure 12) and both negative if $a$ is positive (Figure 13). For both $\mu_{1,2}$ are real. The cubic can have either three real roots (for $c_2 \le c \le c_1$) or only one real root.
\item[{\bf (6)}] $\bm{a^2/4 < b \le a^2/3}$. Thus $\lambda_{1,2}$ are not real and $c_1$ and $c_2$ have the same sign: both are positive if $a < 0$ (Figure 14) and both negative if $a > 0$ (Figure 15). For both $\mu_{1,2}$ are real. The cubic can have either three real roots (for $c_2 \le c \le c_1$) or only one real root.
\item[{\bf (7)}] $\bm{a^2/3 < b}$. Now $\mu_{1,2}$ are not real. In this case, the cubic can have only one real root for all $a$ and all $c$ --- Figure 16 for $a < 0$ and Figure 17 for $a > 0$.
\end{itemize}
\n
Each of these regimes of $b$ is studied for $a < 0$ (even numbered Figures, starting with Figure 4) and $a > 0$ (odd numbered Figures, starting with Figure 5). \\
The depressed cubic (with $a = 0$) is studied separately --- it is on Figures 1, 2, and 3. \\
Note also that the 5 parallel straight lines look like the staves and the intersection points between them and left-hand side $x^2(x+a)$ of (\ref{split}) look like musical notes. A Pythagorean musical analogy (linking proportion to harmony) can be found and the different cubic equations could be endowed with individual tunes --- see \cite{3} for the quartic equation.

\newgeometry{top=40mm}

\section{Algorithm for applying the method}

\begin{itemize}
\item [{\bf (i)}] For the given cubic $x^3 + a x^2 + b x + c$, use $a$ and $b$ to determine which of the above ranges {\bf (1)} to {\bf (7)} for $b$ applies. This will also put the numbers in the list $-c_1$, $-c_0$, $-c_2$, $-ab$, $0$ in increasing order (note that $-c_1 \le c_0 \le c_2$ always).
\item [{\bf (ii)}] Determine the value of $c_0$ from (\ref{c0}), the values of $ab$ and of $\rho_0 = -a/3$, and, using (\ref{rho}), determine the roots $\rho_{1,2}$.
\item [{\bf (iii)}] Solve the three {\it auxiliary quadratic equations} (\ref{q1}), (\ref{q2}), and (\ref{q3}). This will provide $c_{1,2}$, $\mu_{1,2}$, and $\lambda_{1,2}$, respectively. Find also $\xi_{1,2}$ from (\ref{xi}).
\item [{\bf (iv)}] Put the numbers $-c_2$, $-c_1$, $-c_0$, $-ab$, $0$ and the given $-c$ in increasing order.
\item [{\bf (v)}] If $a \ne 0$, skip this point. If the cubic is depressed ($a = 0$) and if $b < 0$, study Figure 1 and, depending on $-c$ in the list from {\bf (i)}, read the isolation intervals of the real roots (which can be either 1 or 3). Figure 2 is a depressed cubic with $b = 0$ and Figure 3 is a depressed cubic with $b > 0$. Both of these cases can have a single root only (except when $b=0$ and $c = 0$ when there is a triple zero root).
\item [{\bf (vi)}] From the sign of the given $a$, determine, depending on the range of $b$, which Figure with number $2n$ applies (should $a < 0$), or which Figure with number $2n+1$ applies (should $a > 0$). Here  $n$ is an integer between $2$ and $8$ inclusive.
\item [{\bf (vii)}] Depending on the position of $-c$ within the list from {\bf (i)}, determine the isolation intervals of the roots of the cubic using the intersection points of  $x^2(x + a)$ and the two parallel lines which are nearest to $- b x - c$. The isolation intervals of the roots are also given in the Figures.
\item [{\bf (viii)}] In the case of $b < 0$ and $c > 0$, the proposed method does not allow one to bind the smallest root $x_3$ from below. This situation is realized on Figure 1 and Figures 4 to 7.  Also, in the case of $b < 0$ and $c < 0$, the biggest root $x_1$ cannot be bound from above. This can be seen on the same Figures. In order to find the isolation intervals of {\it all} roots, in these two cases one needs to find a polynomial root bound: an upper bound $B_U$ and a lower bound $B_L$. This could be any of the many existing root bounds \cite{rs}, \cite{p}, \cite{d}, \cite{k} (which, depending on the coefficients of the polynomial, perform differently). For example, a root bound could be the bigger of 1 and the sum of the absolute values of all coefficients (the Cauchy bound), or one could use the narrower bound \cite{2} which is the bigger of 1 and the sum of the absolute values of all negative coefficients. The bound used in this paper is \cite{d}:  $1 + \sqrt[k]{H}$, where $H$ is the biggest absolute value of all negative coefficients in $x^3 + a x^2 + b x + c$ and $k = 1$ if $a < 0, \,\, k = 2$ if $a > 0$ and $b <0$, and $k = 3$ if $a > 0$ and $b > 0,$ and $c < 0$ (if $a$, $b$, and $c$ are all positive, the upper root bound is zero).
\item [{\bf (ix)}] The ``root harness" (\ref{harness}) provides an additional constraint on the roots: the distance from, say, the biggest root to the smallest root is not smaller than $\sqrt{3}\sqrt{a^2/3 - b}$ and not bigger than $2\sqrt{a^2/3 - b}$. This narrows down the isolation intervals of the roots.
\end{itemize}

\newpage

\newgeometry{top=23mm, bottom=23mm,left=20mm, right=20mm}

\section{Isolation Intervals of the Roots of the General Cubic}

\noindent
The isolation intervals of the roots of the parametric cubic polynomial are explicitly given in the Figures in this Section. All possible cases are presented. \\
Plotted on each Figure is the graph of $x^3 + a x^2$ for $a = 0$ (depressed cubic --- Figures 1 to 3), for $a < 0$ (even numbered Figures from 4 to 16), and for $a > 0$ (odd numbered Figures from 5 to 17), a representative example graph of $- b x - c$, for which $b$ is in its relevant range (depending on $a$) and for some $c$ (note that the free term $c$ can take any value), together with the set of straight lines with slope $-b$, as described in Section 3. These are: $-b x - c_{1,2}$ (when $c_{1,2}$ are real, i.e. for $b \le a^2/3$, the cubic has three real roots if $c_2 \le c \le c_1$), $- b x - c_0, \,\, - b x - a b, \,\,$ and, when relevant, the separatrix $-b x$. \\
To avoid overcrowding, these straight lines are not labelled on the Figures but can be easily identified from their labelled $y$-intercepts. It should also be mentioned that, for ease of reading of the isolation intervals, in all Figures, the points in the $xy$-plane are denoted somewhat unusually. If a point lies on the ordinate, then the name given to the point indicates the value of its $y$-coordinate. For all other points in the plane, the name indicates the $x$-coordinate.\\
\begin{center}
\begin{tabular}{ccc}
\includegraphics[width=45mm]{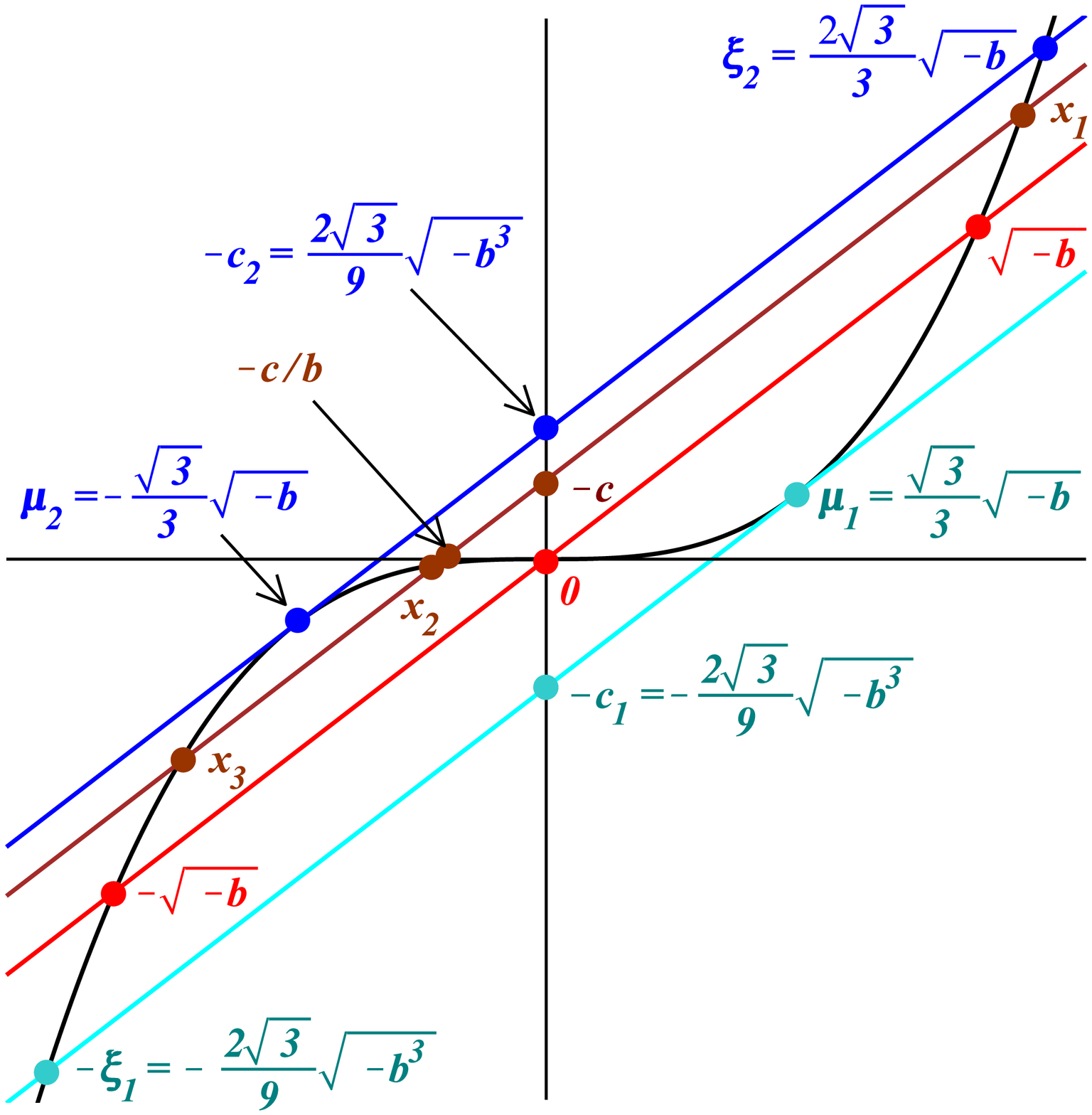} & \includegraphics[width=45mm]{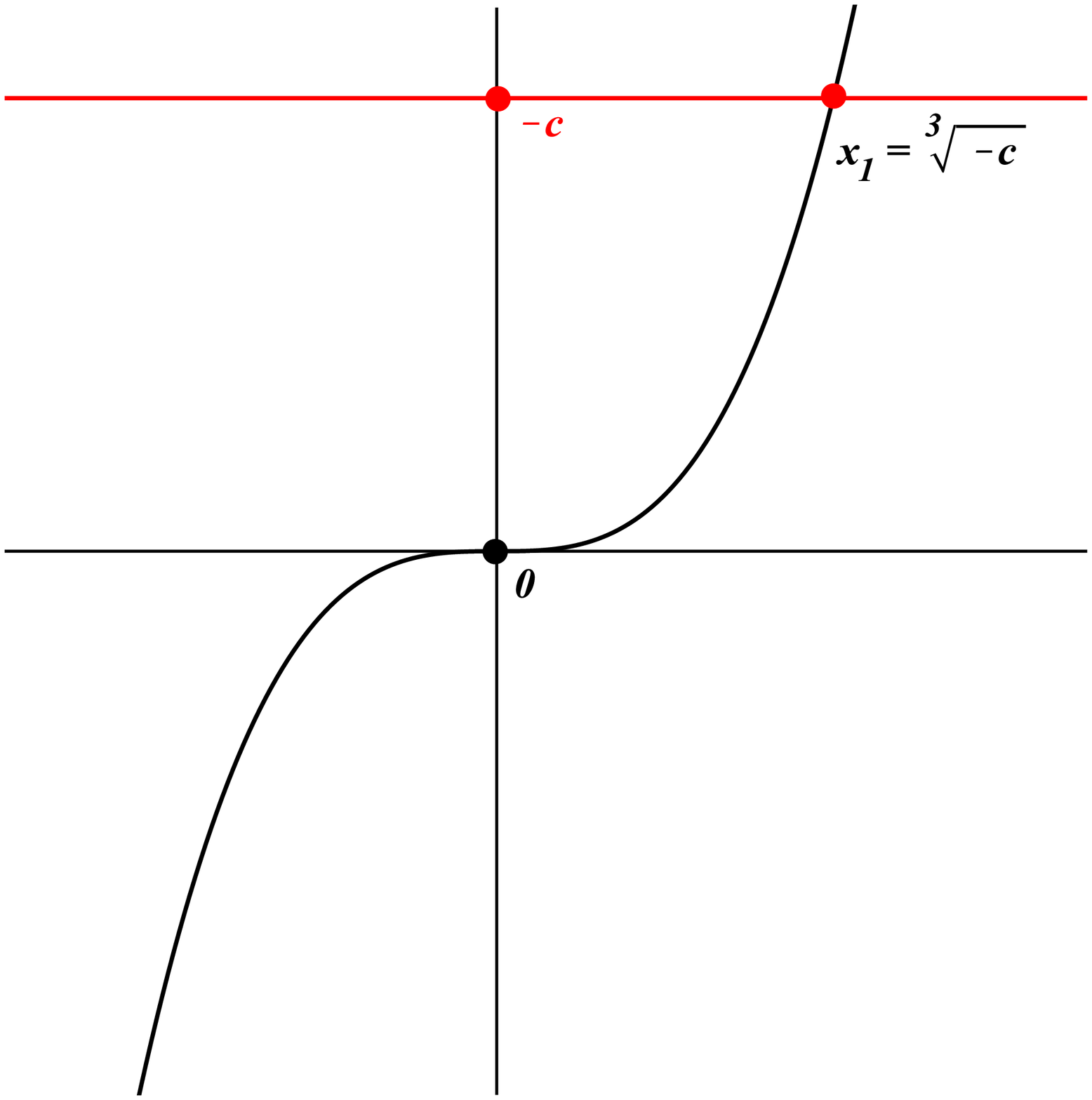} & \includegraphics[width=45mm]{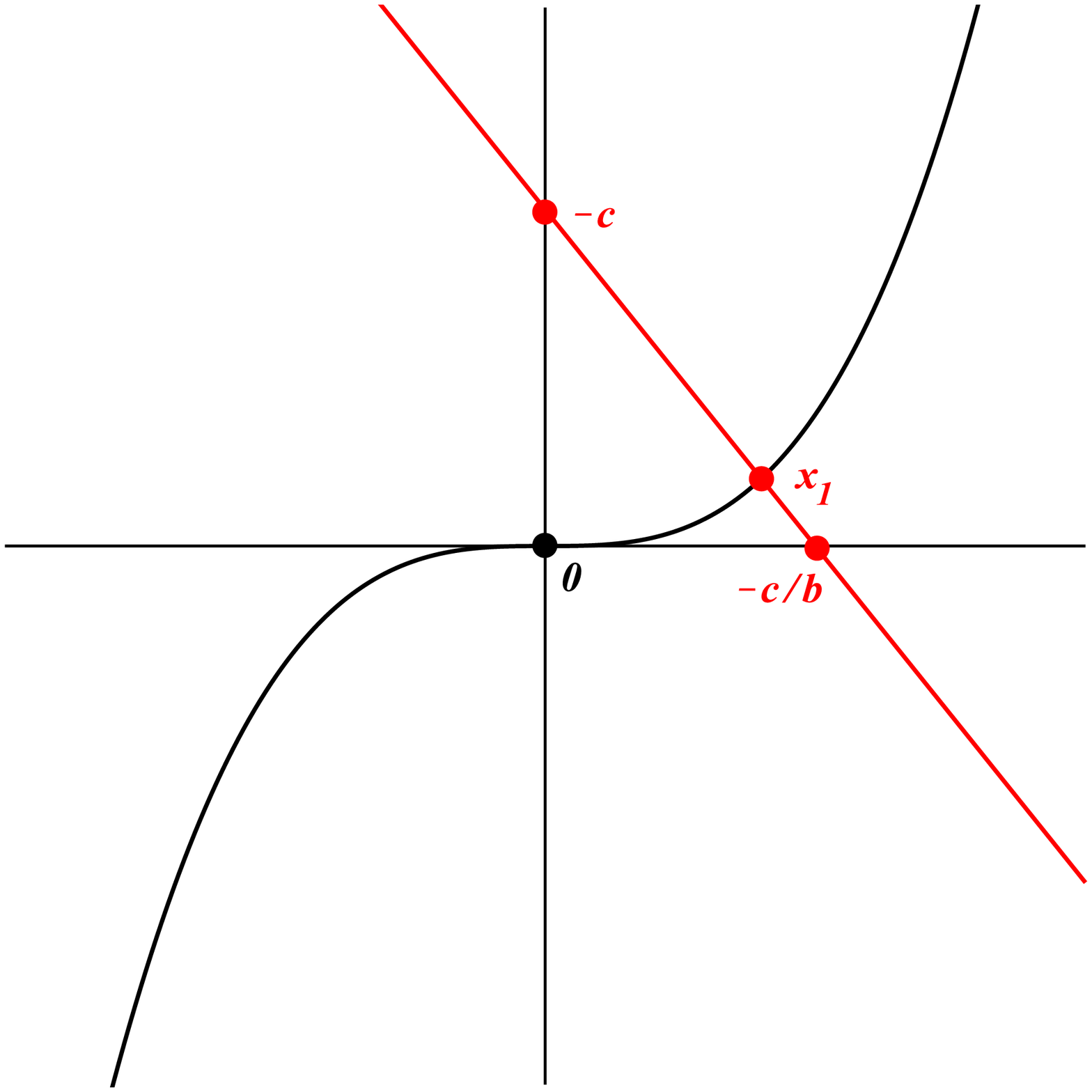} \\
& & \\
{\scriptsize {\bf Figure 1}} & {\scriptsize {\bf Figure 2}} & {\scriptsize {\bf Figure 3}}\\
& & \\
\end{tabular}
\end{center}

\newpage
\newgeometry{top=30mm, bottom=30mm,left=20mm, right=20mm}

\begin{center}
\begin{tabular}{cc}
\includegraphics[width=63mm]{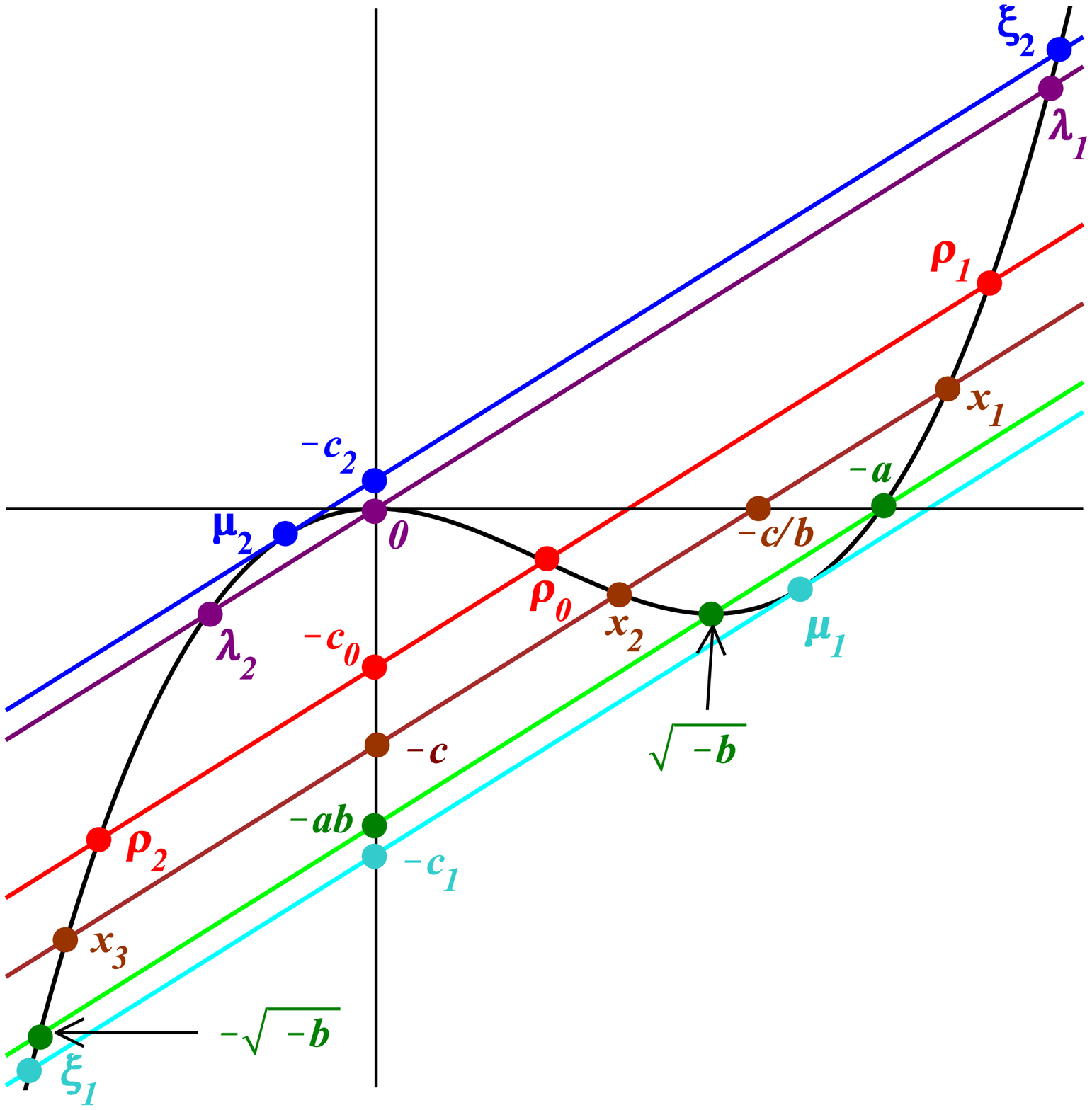} & \includegraphics[width=63mm]{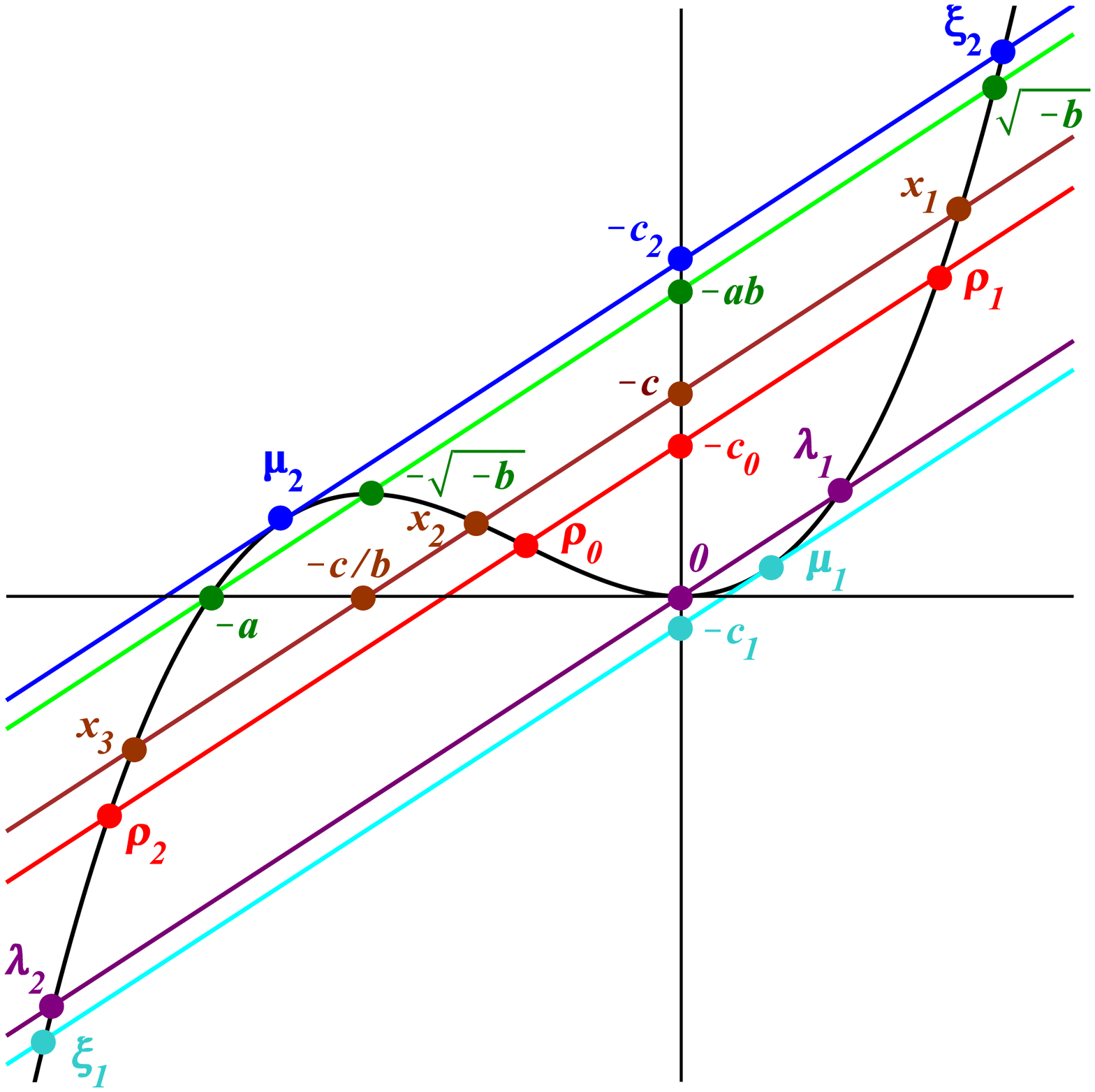} \\
& \\
{\scriptsize {\bf Figure 4}} &  {\scriptsize {\bf Figure 5}} \\
& \\
\end{tabular}
\end{center}

\begin{center}
\begin{tabular}{cc}
\includegraphics[width=63mm]{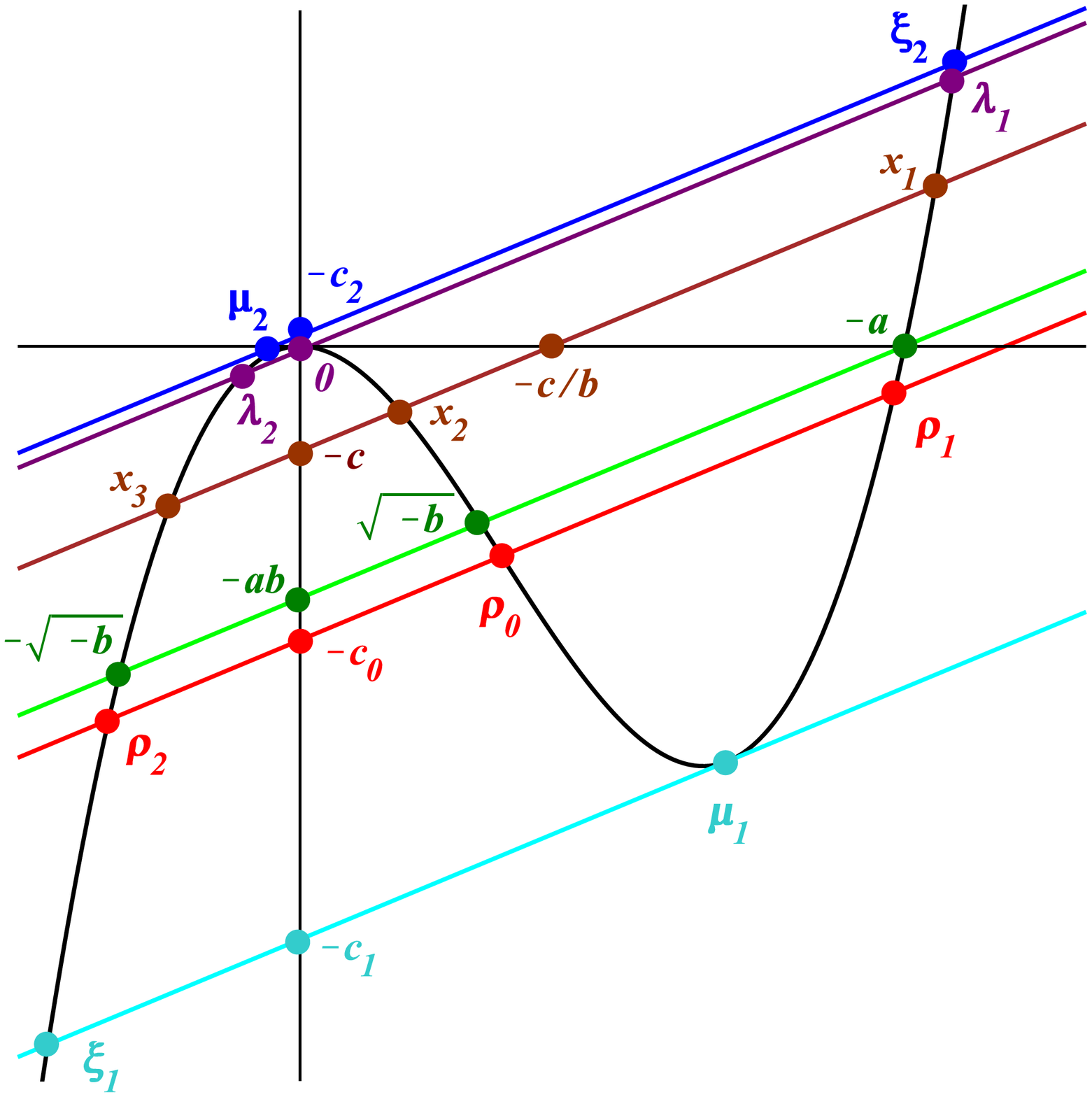} & \includegraphics[width=63mm]{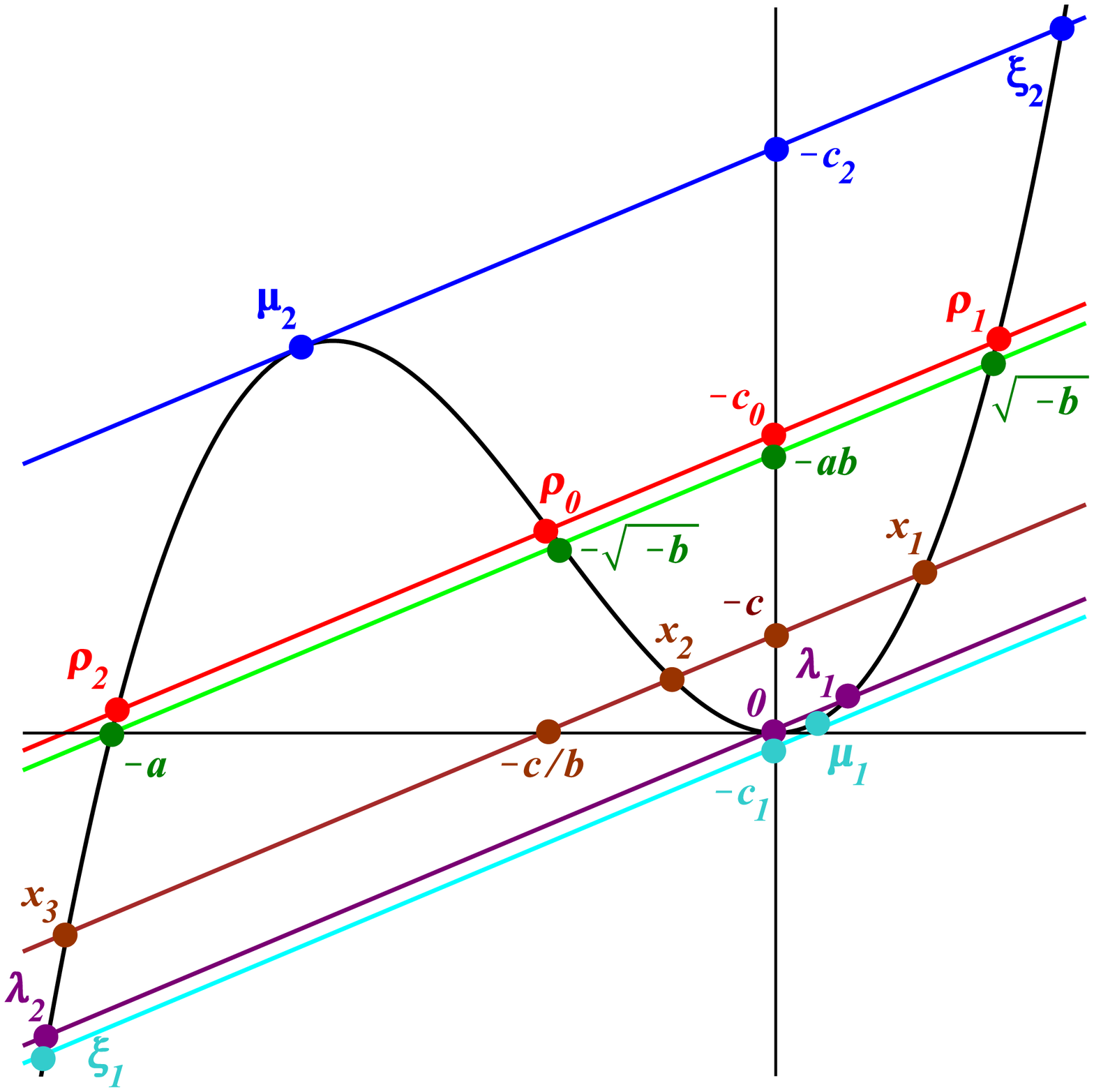} \\
& \\
{\scriptsize {\bf Figure 6}} &  {\scriptsize {\bf Figure 7}} \\
& \\
& \\
\end{tabular}
\end{center}

\begin{center}
\begin{tabular}{cc}
\includegraphics[width=63mm]{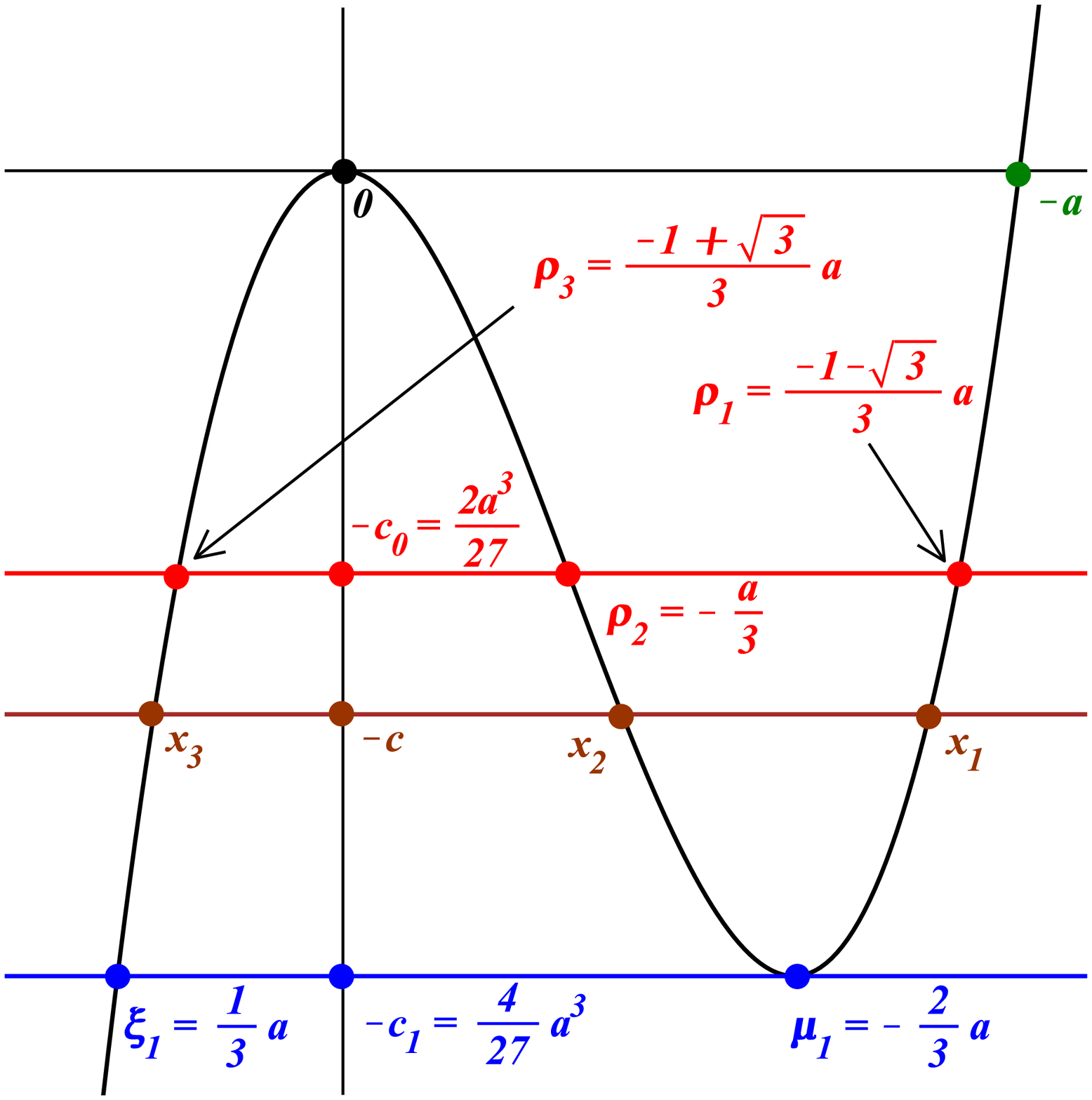} & \includegraphics[width=63mm]{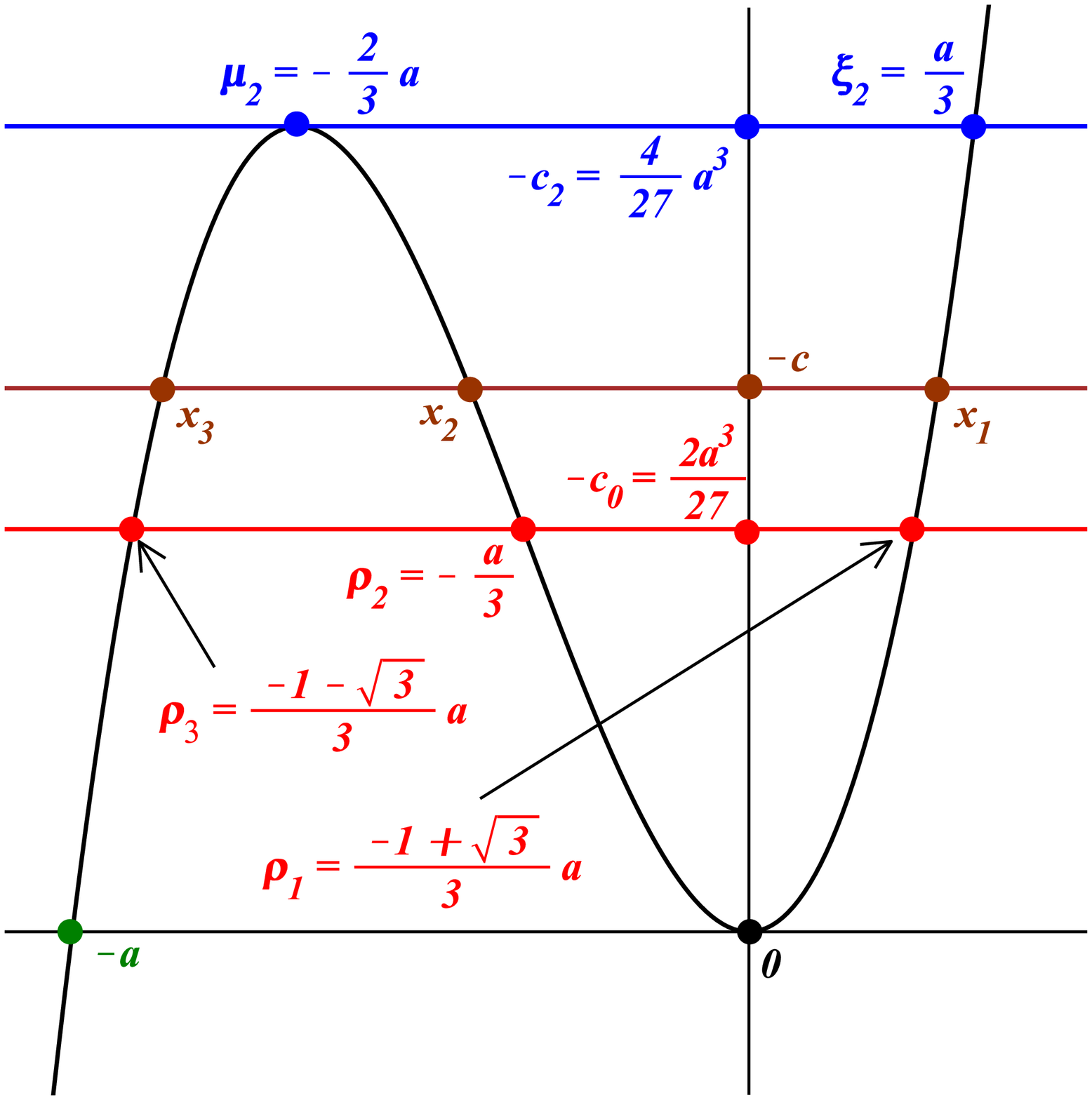} \\
& \\
{\scriptsize {\bf Figure 8}} &  {\scriptsize {\bf Figure 9}} \\
& \\
\end{tabular}
\end{center}

\begin{center}
\begin{tabular}{cc}
\includegraphics[width=63mm]{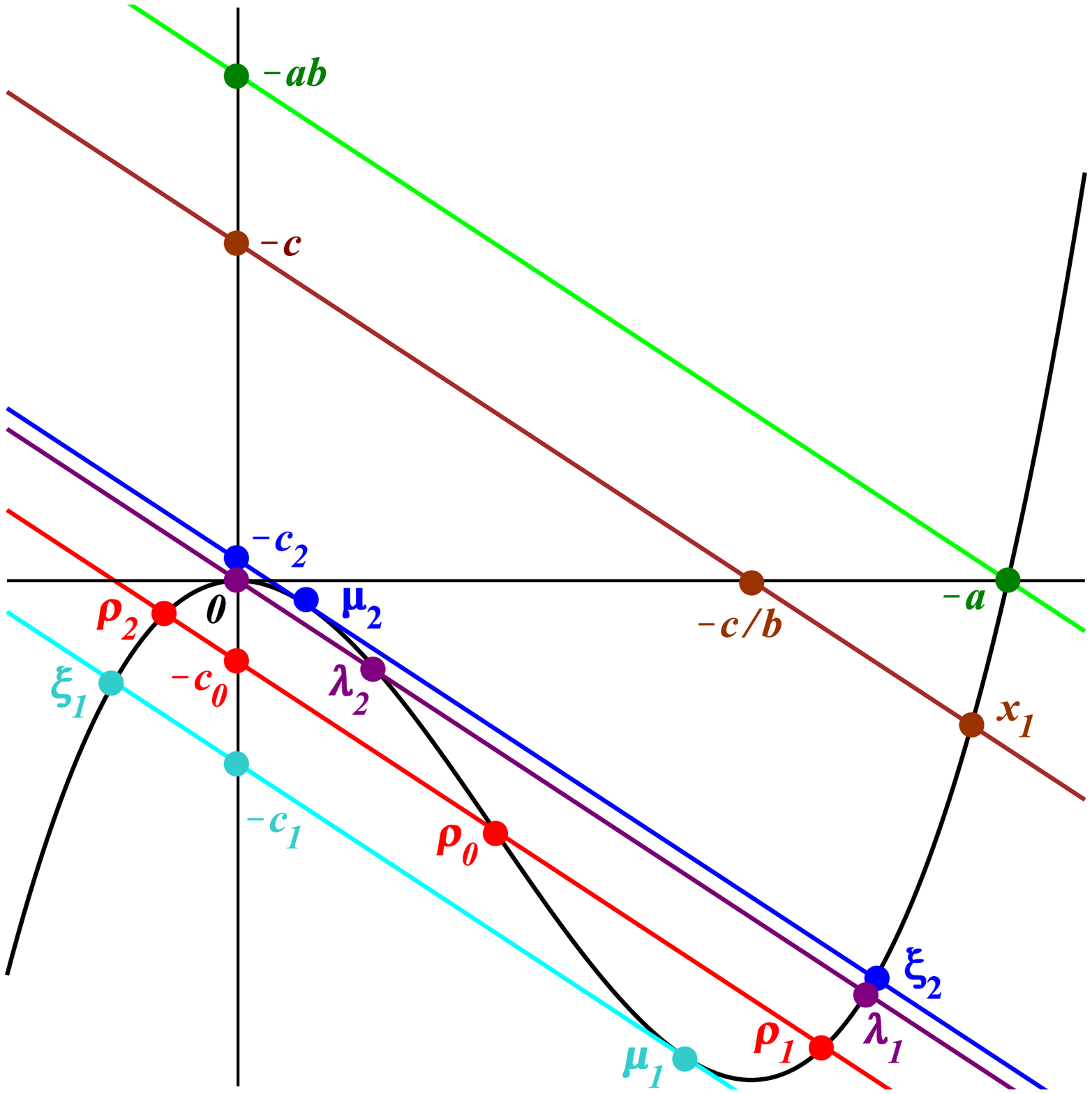} & \includegraphics[width=63mm]{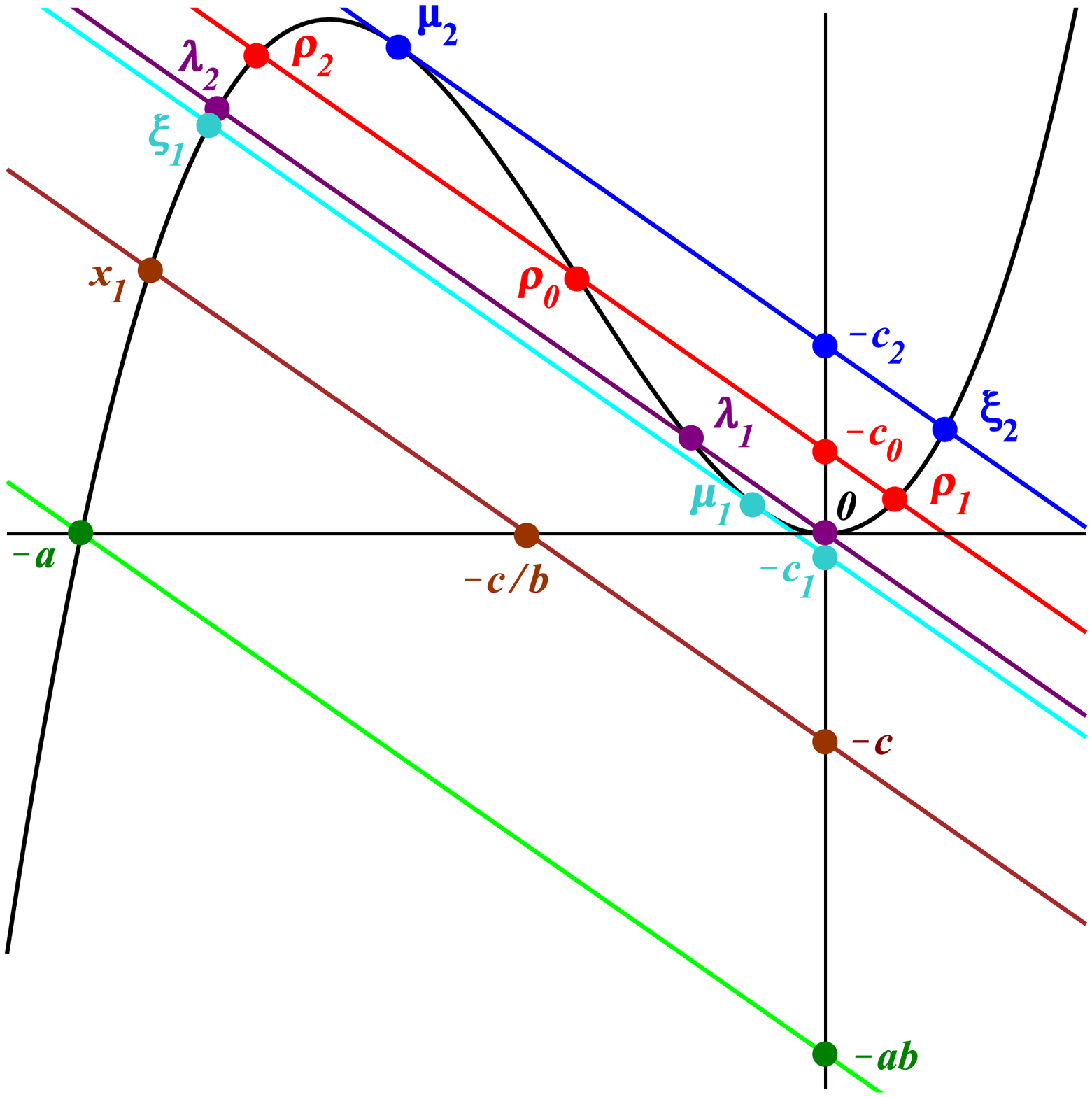} \\
& \\
{\scriptsize {\bf Figure 10}} &  {\scriptsize {\bf Figure 11}} \\
& \\
\end{tabular}
\end{center}

\begin{center}
\begin{tabular}{cc}
\includegraphics[width=63mm]{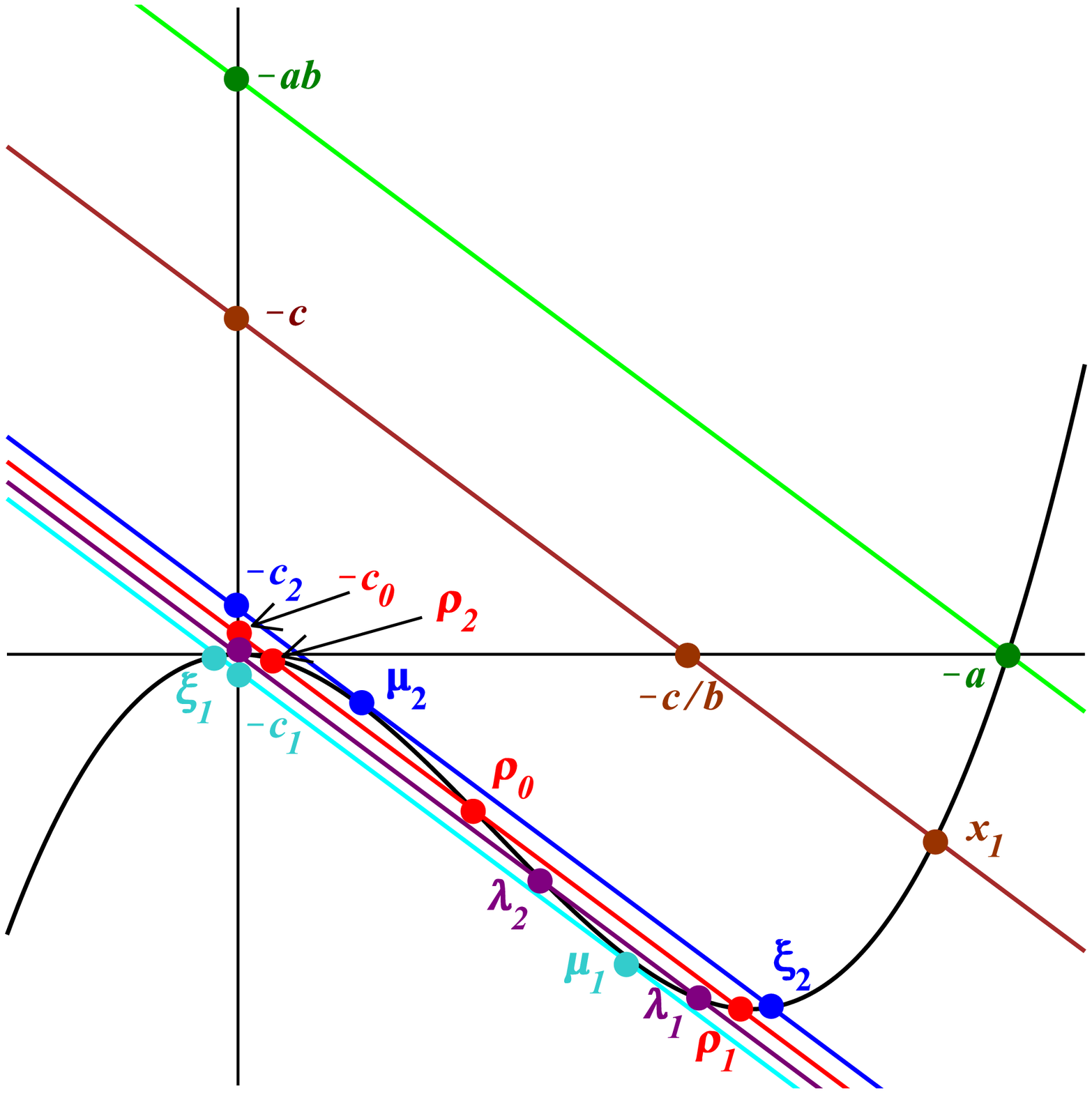} & \includegraphics[width=63mm]{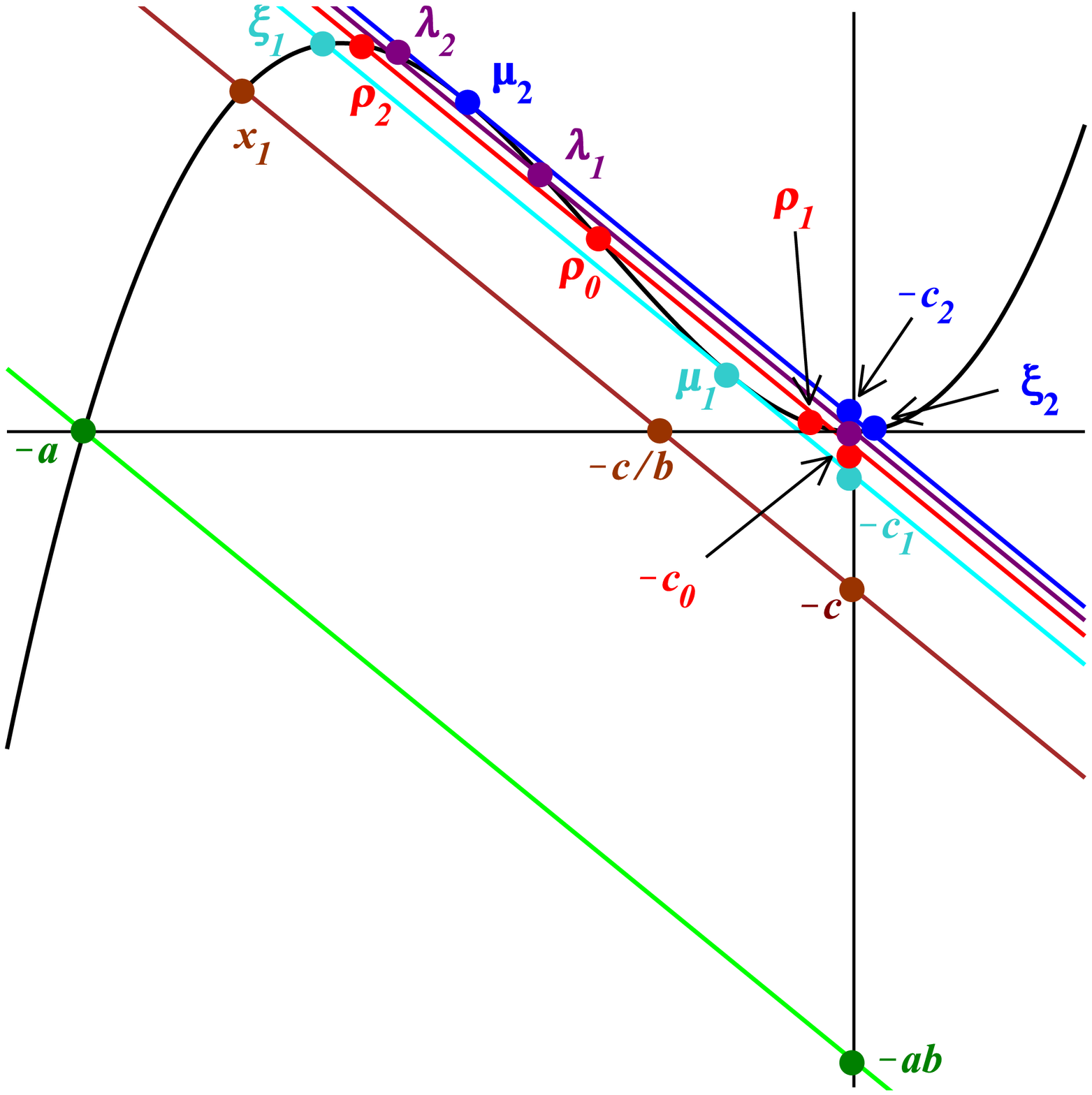} \\
& \\
{\scriptsize {\bf Figure 12}} &  {\scriptsize {\bf Figure 13}} \\
& \\
\end{tabular}
\end{center}

\begin{center}
\begin{tabular}{cc}
\includegraphics[width=63mm]{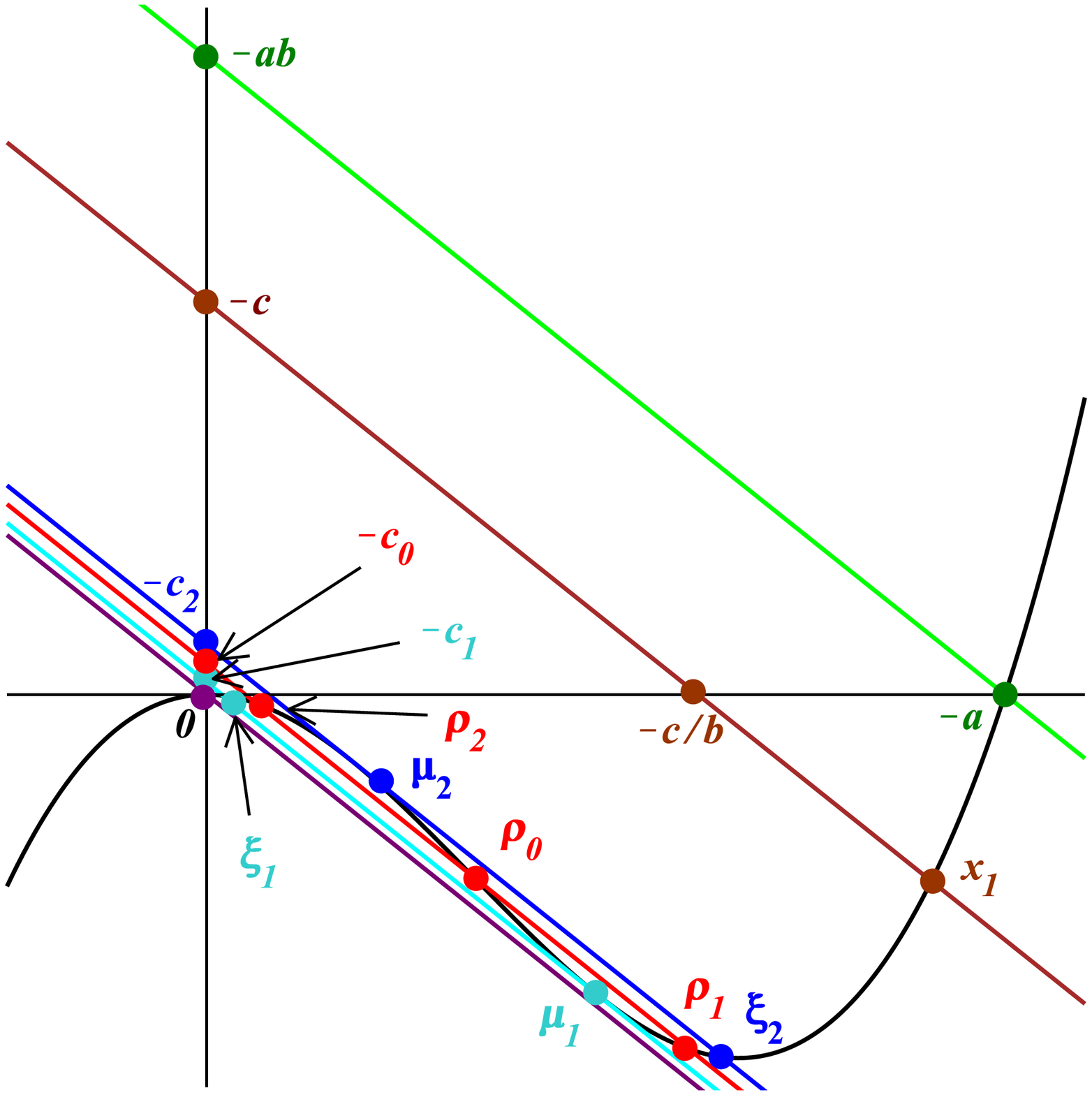} & \includegraphics[width=63mm]{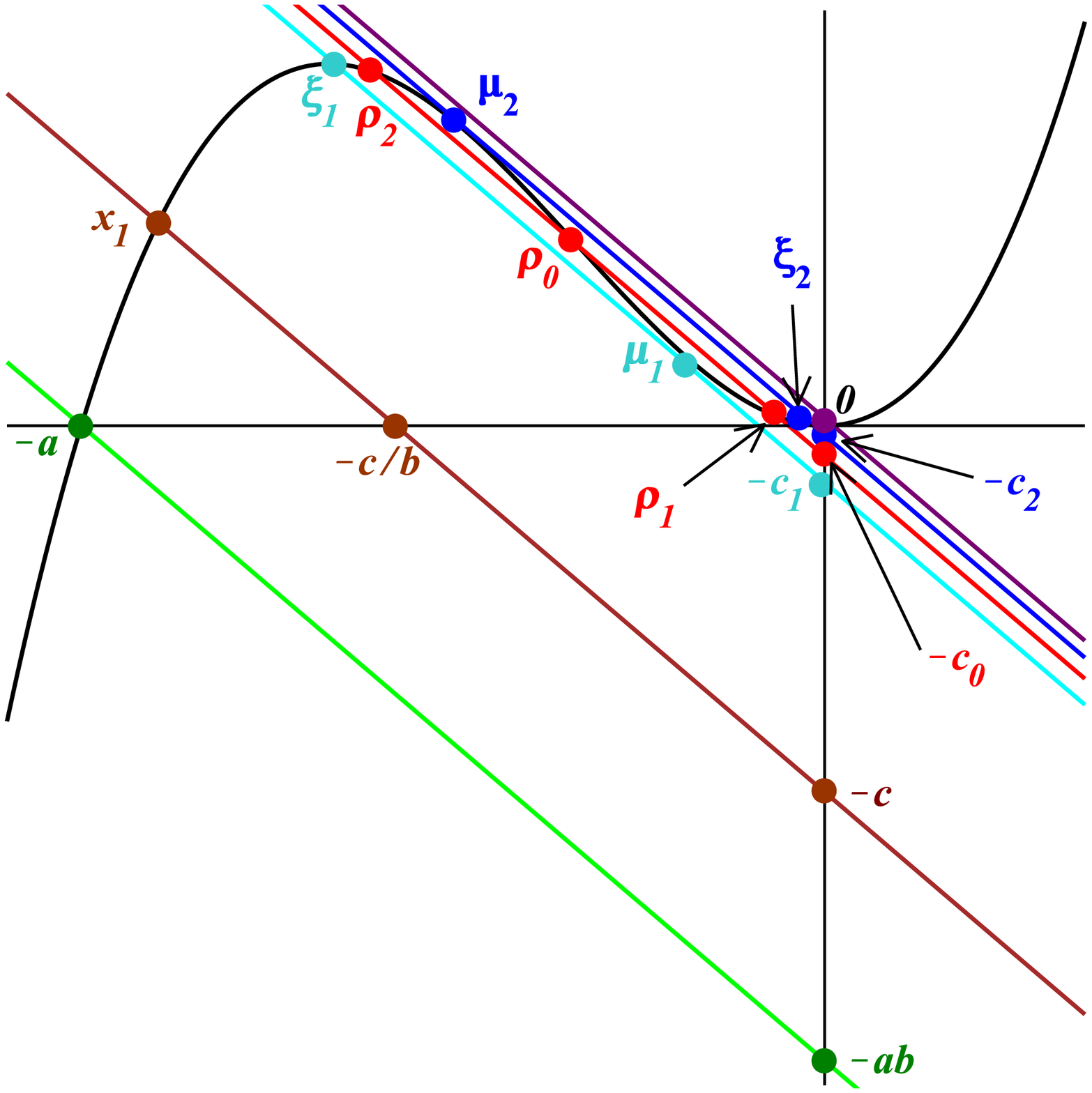} \\
& \\
{\scriptsize {\bf Figure 14}} &  {\scriptsize {\bf Figure 15}} \\
& \\
\end{tabular}
\end{center}

\begin{center}
\begin{tabular}{cc}
\includegraphics[width=63mm]{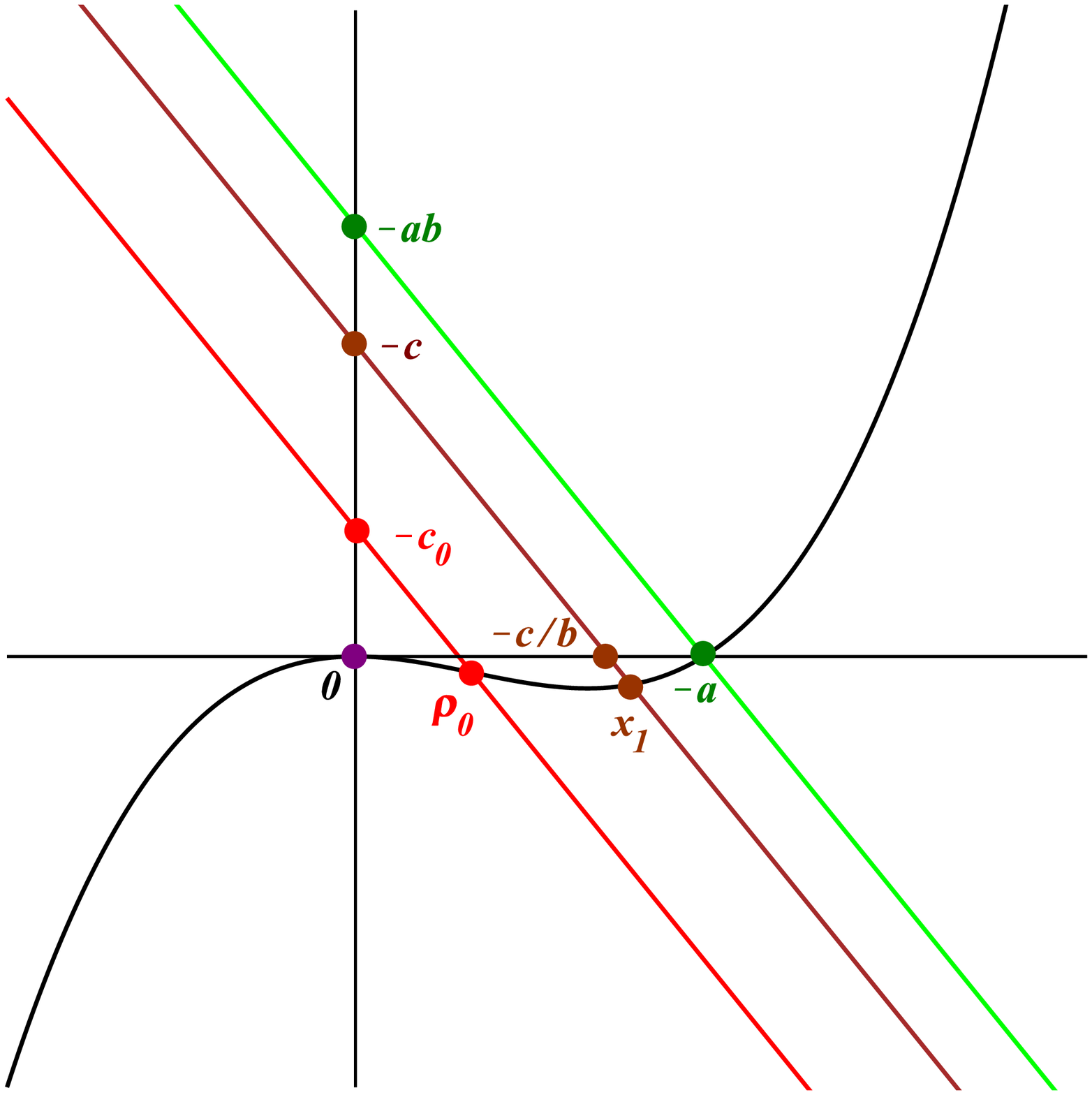} & \includegraphics[width=63mm]{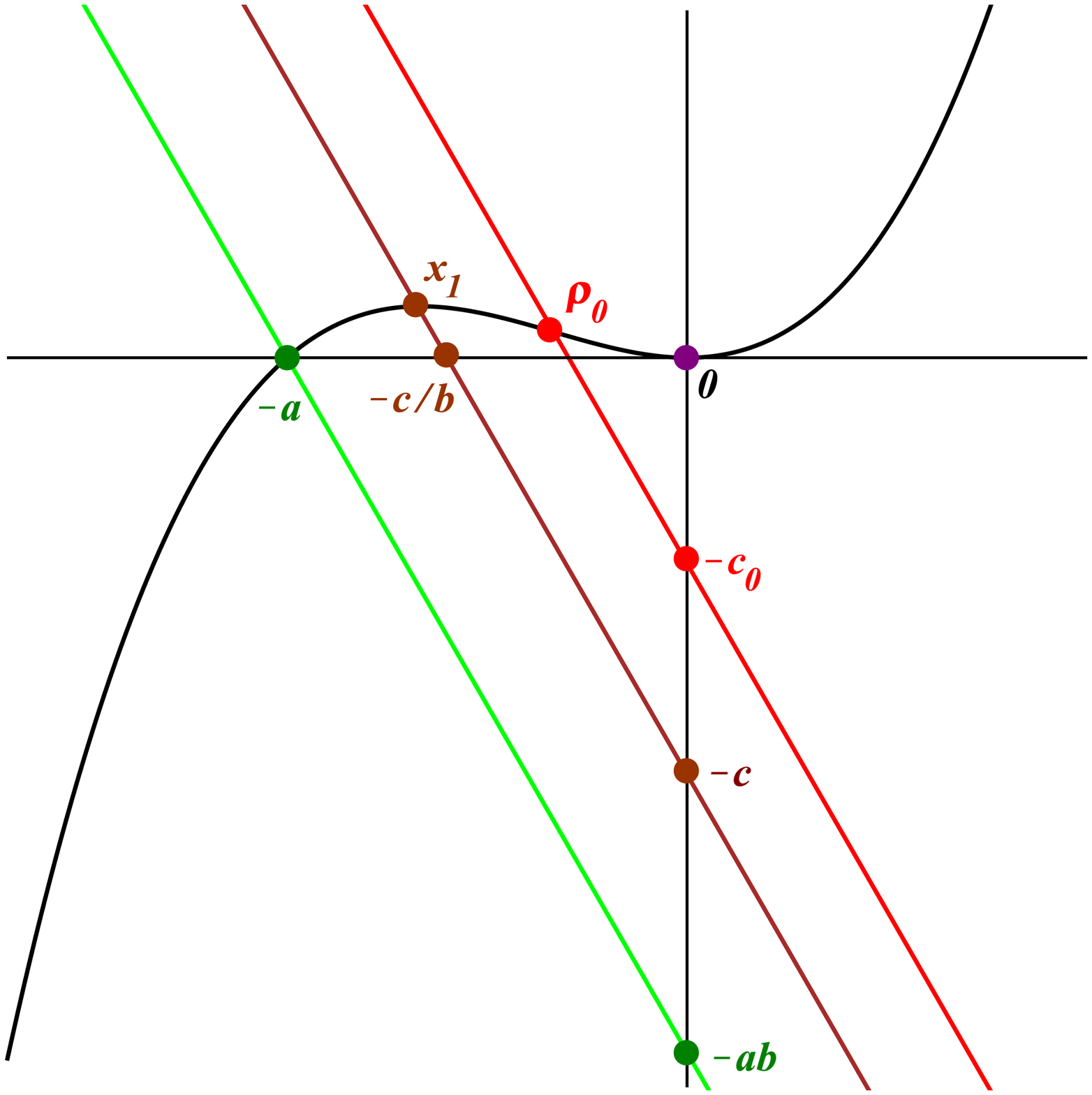} \\
& \\
{\scriptsize {\bf Figure 16}} &  {\scriptsize {\bf Figure 17}} \\
& \\
\multicolumn{1}{c}{\begin{minipage}{17em}
\scriptsize
\begin{center}
\vskip-.3cm
\underline{$\bm{b > a^2/3, \quad a < 0}$}
\end{center}
\begin{itemize}
\item[{\bf (1)}] $\bm{- c < 0}$ --- one negative root:
\begin{itemize}
\item[{\bf (i)}] $-c/b < x_1 < 0$
\end{itemize}
\item[{\bf (2)}] $\bm{0 \le - c < - c_0}$ --- one non-negative root:
\begin{itemize}
\item[{\bf (i)}] $-c/b \le x_1 <  \rho_0$
\end{itemize}
\item[{\bf (3)}] $\bm{ - c_0  \le - c < - ab}$ (shown on graph) --- one positive root:
\begin{itemize}
\item[{\bf (i)}] max$\{ -c/b, \rho_0 \}  \le  x_1 <  -a$
\end{itemize}
\item[{\bf (4)}] $\bm{- ab \le - c }$ --- one positive root:
\begin{itemize}
\item[{\bf (i)}] $-a \le x_1 < -c/b$
\end{itemize}
\end{itemize}
\end{minipage}}
& \multicolumn{1}{c}{\begin{minipage}{17em}
\scriptsize
\begin{center}
\vskip-.3cm
\underline{$\bm{b > a^2/3, \quad a > 0}$}
\end{center}
\begin{itemize}
\item[{\bf (1)}] $\bm{- c < -ab}$ --- one negative root:
\begin{itemize}
\item[{\bf (i)}] $-c/b < x_1 < -a$
\end{itemize}
\item[{\bf (2)}] $\bm{- ab \le - c < - c_0}$ (shown on graph)  --- one negative root:
\begin{itemize}
\item[{\bf (i)}] $-a \le x_1 < \,\, $min$\{ -c/b, \rho_0 \}$
\end{itemize}
\item[{\bf (3)}] $\bm{ - c_0  \le - c < 0}$ --- one negative root:
\begin{itemize}
\item[{\bf (i)}] $\rho_0  \le  x_1 <  - c/b$
\end{itemize}
\item[{\bf (4)}] $\bm{0 \le - c}$ --- one non-negative root:
\begin{itemize}
\item[{\bf (i)}] $0 \le x_1 < -c/b$
\end{itemize}
\end{itemize}
\end{minipage}}
\\
\end{tabular}
\end{center}

\newpage

\newgeometry{top=35mm, bottom=35mm,left=23mm, right=23mm}

\section{Improved Complete Root Classification for the General Cubic Polynomial}

\n
Excluding the case of a zero root (which occurs when $c = 0$ and will not be considered), all possible cases of the {\it Improved Complete Root Classification} of the cubic polynomial, with their conditions, relevant Figures, and corresponding root isolation intervals, are presented in the following six tables:

\vskip.5cm

\begin{table}[!h]
\hskip.95cm
\begin{tabular}{|ccccc|}
\hline
\multicolumn{5}{|c|}{{\bf (I) Three Positive Roots} $\bm{(0 < x_3 \le x_2 \le x_1)}$} \\ \hline
\multicolumn{1}{|c|}{$\bm{b}$} & \multicolumn{1}{c|}{$\bm{a}$} & \multicolumn{1}{c|}{$\bm{c}$} & \multicolumn{1}{c|}{{\bf Figure}} & {\bf Root Isolation Intervals} \\ \hline
\multicolumn{1}{|c|}{\multirow{3}{*}{$0 < b \le 2a^2/9$}} & \multicolumn{1}{c|}{\multirow{3}{*}{$a < 0$}} & \multicolumn{1}{c|}{\multirow{3}{*}{$0 < -c \le -c_2$}} & \multicolumn{1}{c|}{\multirow{3}{*}{10 (4)}} &  $-c/b \le x_3 \le \mu_2,$ \\
\multicolumn{1}{|c|}{} & \multicolumn{1}{c|}{} & \multicolumn{1}{c|}{} & \multicolumn{1}{c|}{}  & $\mu_2 \le x_2 \le \lambda_2$,  \\
\multicolumn{1}{|c|}{} & \multicolumn{1}{c|}{} & \multicolumn{1}{c|}{} & \multicolumn{1}{c|}{}  & $\lambda_1 \le x_1 \le \xi_2$. \\ \hline
\multicolumn{1}{|c|}{\multirow{6}{*}{$2a^2/9 < b \le a^2/4$}} & \multicolumn{1}{c|}{\multirow{6}{*}{$a < 0$}} & \multicolumn{1}{c|}{\multirow{3}{*}{$0 < -c < -c_0$}} & \multicolumn{1}{c|}{\multirow{3}{*}{12 (3)}} &  $-c/b \le x_3 < \rho_2,$ \\
\multicolumn{1}{|c|}{} & \multicolumn{1}{c|}{} & \multicolumn{1}{c|}{} & \multicolumn{1}{c|}{} & $\rho_0 < x_2 \le \lambda_2,$ \\
\multicolumn{1}{|c|}{} & \multicolumn{1}{c|}{} & \multicolumn{1}{c|}{} & \multicolumn{1}{c|}{} & $\lambda_1 \le x_1 < \rho_1.$ \\
\hline
\end{tabular}
\end{table}

\vskip.5cm

\begin{table}[!h]
\hskip.95cm
\begin{tabular}{|ccccc|}
\hline
\multicolumn{5}{|c|}{{\bf (II) Three Negative Roots} $\bm{(x_3 \le x_2 \le x_1 < 0)}$} \\ \hline
\multicolumn{1}{|c|}{$\bm{b}$} & \multicolumn{1}{c|}{$\bm{a}$} & \multicolumn{1}{c|}{$\bm{c}$} & \multicolumn{1}{c|}{{\bf Figure}} & {\bf Root Isolation Intervals} \\ \hline
\multicolumn{1}{|c|}{\multirow{3}{*}{$0 < b \le 2a^2/9$}} & \multicolumn{1}{c|}{\multirow{3}{*}{$a > 0$}} & \multicolumn{1}{c|}{\multirow{3}{*}{$- c_1 \le - c < 0$}} & \multicolumn{1}{c|}{\multirow{3}{*}{11 (3)}} &  $\xi_1  \le  x_3 < \lambda_2,$ \\
\multicolumn{1}{|c|}{} & \multicolumn{1}{c|}{} & \multicolumn{1}{c|}{} & \multicolumn{1}{c|}{}  & $\lambda_1 < x_2 \le \mu_1$,  \\
\multicolumn{1}{|c|}{} & \multicolumn{1}{c|}{} & \multicolumn{1}{c|}{} & \multicolumn{1}{c|}{}  & $\mu_1  \le x_1 < - c/b$. \\ \hline
\multicolumn{1}{|c|}{\multirow{6}{*}{$2a^2/9 < b \le a^2/4$}} & \multicolumn{1}{c|}{\multirow{6}{*}{$a > 0$}} & \multicolumn{1}{c|}{\multirow{3}{*}{$-c_1 \le -c < -c_0$}} & \multicolumn{1}{c|}{\multirow{3}{*}{13 (3)}} &  $\xi_1 \le x_3 < \rho_2,$ \\
\multicolumn{1}{|c|}{} & \multicolumn{1}{c|}{} & \multicolumn{1}{c|}{} & \multicolumn{1}{c|}{} & $\rho_0 < x_2 \le \mu_1,$ \\
\multicolumn{1}{|c|}{} & \multicolumn{1}{c|}{} & \multicolumn{1}{c|}{} & \multicolumn{1}{c|}{} & $\mu_1 \le x_1 < $ min$\{ -c/b, \rho_1 \}.$ \\
\hline
\end{tabular}
\end{table}

\newpage

\newgeometry{left=23mm, right=23mm}

\begin{table}[!h]
\begin{tabular}{|ccccc|}
\hline
\multicolumn{5}{|c|}{{\bf (III) Two Positive Roots} $(\bm{0 < x_2 \le x_1})$ {\bf and One Negative Root} $(\bm{x_3 < 0})$} \\ \hline
\multicolumn{1}{|c|}{$\bm{b}$} & \multicolumn{1}{c|}{$\bm{a}$} & \multicolumn{1}{c|}{$\bm{c}$} & \multicolumn{1}{c|}{{\bf Figure}} & {\bf Root Isolation Intervals} \\ \hline
\multicolumn{1}{|c|}{\multirow{3}{*}{$b < 0$}} & \multicolumn{1}{c|}{\multirow{3}{*}{$a = 0$}} & \multicolumn{1}{c|}{\multirow{3}{*}{$- (2\sqrt{3}/9)\sqrt{-b^3} \le -c < 0$}} & \multicolumn{1}{c|}{\multirow{3}{*}{1 (2)}} &  $-(2\sqrt{3}/3)\sqrt{-b} \le x_3 < -\sqrt{-b},$ \\
\multicolumn{1}{|c|}{} & \multicolumn{1}{c|}{} & \multicolumn{1}{c|}{} & \multicolumn{1}{c|}{}  & $-c/b < x_2 \le (\sqrt{3}/3)\sqrt{-b}$,  \\
\multicolumn{1}{|c|}{} & \multicolumn{1}{c|}{} & \multicolumn{1}{c|}{} & \multicolumn{1}{c|}{}  & $(\sqrt{3}/3)\sqrt{-b} \le x_1 < \sqrt{-b}$. \\ \hline
\multicolumn{1}{|c|}{\multirow{9}{*}{$b < -a^2/9$}} & \multicolumn{1}{c|}{\multirow{9}{*}{$a < 0$}} & \multicolumn{1}{c|}{\multirow{3}{*}{$-c_1 \le -c < - ab$}} & \multicolumn{1}{c|}{\multirow{3}{*}{4 (2)}} &  $\xi_1 \le x_3 < -\sqrt{-b},$ \\
\multicolumn{1}{|c|}{} & \multicolumn{1}{c|}{} & \multicolumn{1}{c|}{} & \multicolumn{1}{c|}{} & $\sqrt{-b} < x_2 \le \mu_1,$ \\
\multicolumn{1}{|c|}{} & \multicolumn{1}{c|}{} & \multicolumn{1}{c|}{} & \multicolumn{1}{c|}{} & $\mu_1 \le x_1 < -a.$ \\
\hline
\end{tabular}
\end{table}

\newpage

\newgeometry{left=23mm, right=23mm}

\begin{table}[!h]
\begin{tabular}{|ccccc|}
\hline
\multicolumn{5}{|c|}{{\bf (IV) One Positive Root} $(\bm{x_1 > 0})$  {\bf and Two Negative Roots} $(\bm{x_3 \le x_2 < 0})$} \\ \hline
\multicolumn{1}{|c|}{$\bm{b}$} & \multicolumn{1}{c|}{$\bm{a}$} & \multicolumn{1}{c|}{$\bm{c}$} & \multicolumn{1}{c|}{{\bf Figure}} & {\bf Root Isolation Intervals} \\ \hline
\multicolumn{1}{|c|}{\multirow{3}{*}{$b < 0$}} & \multicolumn{1}{c|}{\multirow{3}{*}{$a = 0$}} & \multicolumn{1}{c|}{\multirow{3}{*}{$0 < -c \le (2\sqrt{3}/9)\sqrt{-b^3}$}} & \multicolumn{1}{c|}{\multirow{3}{*}{1 (3)}} &  $-\sqrt{-b} \le x_3 \le - (\sqrt{3}/3)\sqrt{-b},$ \\
\multicolumn{1}{|c|}{} & \multicolumn{1}{c|}{} & \multicolumn{1}{c|}{} & \multicolumn{1}{c|}{}  & $- (\sqrt{3}/3)\sqrt{-b} \le x_2 \le -c/b$,  \\
\multicolumn{1}{|c|}{} & \multicolumn{1}{c|}{} & \multicolumn{1}{c|}{} & \multicolumn{1}{c|}{}  & $\sqrt{-b} \le x_1 < (2\sqrt{3}/3)\sqrt{-b}$. \\ \hline
\multicolumn{1}{|c|}{\multirow{3}{*}{$b < - a^2/9$}} & \multicolumn{1}{c|}{\multirow{3}{*}{$a < 0$}} & \multicolumn{1}{c|}{\multirow{3}{*}{$0 < -c \le -c_2$}} & \multicolumn{1}{c|}{\multirow{3}{*}{4 (5)}} &  $\lambda_2 \le x_3 \le \mu_2,$ \\
\multicolumn{1}{|c|}{} & \multicolumn{1}{c|}{} & \multicolumn{1}{c|}{} & \multicolumn{1}{c|}{}  & $\mu_2 \le x_2 \le -c/b$,  \\
\multicolumn{1}{|c|}{} & \multicolumn{1}{c|}{} & \multicolumn{1}{c|}{} & \multicolumn{1}{c|}{}  & $\lambda_1 \le x_1 \le \xi_2$. \\
\hline
\end{tabular}
\end{table}

\newpage

\newgeometry{top=23mm, bottom=23mm,left=23mm, right=23mm}

\begin{table}[!h]
\begin{tabular}{|ccccc|}
\hline
\multicolumn{5}{|c|}{{\bf (V) One Positive Root} $(\bm{x_1 > 0})$ {\bf and a Pair of Complex-Conjugate Roots}} \\ \hline
\multicolumn{1}{|c|}{$\bm{b}$} & \multicolumn{1}{c|}{$\bm{a}$} & \multicolumn{1}{c|}{$\bm{c}$} & \multicolumn{1}{c|}{{\bf Figure}} & {\bf Root Isolation Intervals} \\ \hline
\multicolumn{1}{|c|}{$b < 0$} & \multicolumn{1}{c|}{$a = 0$} & \multicolumn{1}{c|}{$(2\sqrt{3}/9)\sqrt{-b^3} < -c$} & \multicolumn{1}{c|}{1 (4)} & $(2 \sqrt{3}/3) \sqrt{-b}  < x_1  < 1 \, + \,$max$\{|b|, |c| \}$ \\ \hline
\multicolumn{1}{|c|}{$b = 0$} & \multicolumn{1}{c|}{$a = 0$} & \multicolumn{1}{c|}{$-c > 0$} & \multicolumn{1}{c|}{2 (3)} & $x_1 = \sqrt[3]{-c}$ \\ \hline
\multicolumn{1}{|c|}{$b > 0$} & \multicolumn{1}{c|}{$a = 0$} & \multicolumn{1}{c|}{$-c > 0$} & \multicolumn{1}{c|}{3 (2)} & $0 <  x_1  < -c/b$ \\ \hline
\multicolumn{1}{|c|}{$b < - a^2/9$} & \multicolumn{1}{c|}{$a < 0$} & \multicolumn{1}{c|}{$-c_2 < -c$} & \multicolumn{1}{c|}{4 (6)} & $\xi_2 < x_1 < 1 \, + \, $max$\{|a|, |b|, |c| \}$ \\ \hline
\multicolumn{1}{|c|}{$b < - a^2/9$} & \multicolumn{1}{c|}{$a > 0$} & \multicolumn{1}{c|}{$-c_2 < -c$} & \multicolumn{1}{c|}{5 (6)} & $\xi_2 < x_1 < 1 \, + \, \sqrt{\mbox{max}\{|b|, |c| \}}$ \\ \hline
\multicolumn{1}{|c|}{$- a^2/9 \le b < 0$} & \multicolumn{1}{c|}{$a < 0$} & \multicolumn{1}{c|}{$-c_2 < -c$} & \multicolumn{1}{c|}{6 (6)} & $\xi_2 < x_1 < 1 \, + \, $max$\{|a|, |b|, |c| \}$ \\
\hline
\end{tabular}
\end{table}

\begin{table}[!h]
\hskip-.4cm
\begin{tabular}{|ccccc|}
\hline
\multicolumn{5}{|c|}{{\bf (VI) One Negative Root} $(\bm{x_1 < 0})$ {\bf and a Pair of Complex-Conjugate Roots}} \\ \hline
\multicolumn{1}{|c|}{$\bm{b}$} & \multicolumn{1}{c|}{$\bm{a}$} & \multicolumn{1}{c|}{$\bm{c}$} & \multicolumn{1}{c|}{{\bf Figure}} & {\bf Root Isolation Intervals} \\ \hline
\multicolumn{1}{|c|}{$b < 0$} & \multicolumn{1}{c|}{$a = 0$} & \multicolumn{1}{c|}{$-c < -(2\sqrt{3}/9)\sqrt{-b^3}$} & \multicolumn{1}{c|}{1 (1)} & $-(1 \, + \,$max$\{|b|, c \}) < x_1 < -(2 \sqrt{3}/3) \sqrt{-b}$ \\ \hline
\multicolumn{1}{|c|}{$b = 0$} & \multicolumn{1}{c|}{$a = 0$} & \multicolumn{1}{c|}{$-c < 0$} & \multicolumn{1}{c|}{2 (1)} & $x_1 = \sqrt[3]{-c}$ \\ \hline
\multicolumn{1}{|c|}{$b > 0$} & \multicolumn{1}{c|}{$a = 0$} & \multicolumn{1}{c|}{$-c < 0$} & \multicolumn{1}{c|}{3 (1)} & $-c/b <  x_1  < 0$ \\ \hline
\multicolumn{1}{|c|}{$b < - a^2/9$} & \multicolumn{1}{c|}{$a < 0$} & \multicolumn{1}{c|}{$-c < -c_1$} & \multicolumn{1}{c|}{4 (1)} & $-(1 \, + \, \sqrt{\mbox{max} \{|b|, c \}}) < x_1 < \xi_1$ \\ \hline
\multicolumn{1}{|c|}{$b < - a^2/9$} & \multicolumn{1}{c|}{$a > 0$} & \multicolumn{1}{c|}{$-c < -c_1$} & \multicolumn{1}{c|}{5 (1)} & $-(1 \, + \,$max$\{a, |b|, c \}) < x_1 < \xi_1$ \\ \hline
\multicolumn{1}{|c|}{$- a^2/9 \le b < 0$} & \multicolumn{1}{c|}{$a < 0$} & \multicolumn{1}{c|}{$-c < -c_1$} & \multicolumn{1}{c|}{6 (1)} & $-(1 \, + \, \sqrt{\mbox{max} \{|b|, c \}}) < x_1 < \xi_1$ \\
\hline
\end{tabular}
\end{table}

\newpage

\restoregeometry

\noindent
If one is not interested in the root isolation intervals, recalling that
\b
c_{1,2}(a,b) =  -\frac{2}{27} \, a^3 + \frac{1}{3} \, a b \, \pm \, \frac{2}{27} \, \sqrt{( a^2 - 3 b)^3}
\e
(as only those are needed), the {\it improved complete root classification} can be summarized as follows:

\vskip.5cm

\begin{table}[!h]
\begin{center}
\begin{tabular}{|ccc|}
\hline
\multicolumn{3}{|c|}{{\bf (I) Three Positive Roots}}  \\ \hline
\multicolumn{1}{|c|}{$\bm{a}$} & \multicolumn{1}{c|}{$\bm{b}$} & $\bm{c}$ \\ \hline
\multicolumn{1}{|c|}{\multirow{2}{*}{$a < 0$}} & \multicolumn{1}{c|}{$0 < b \le a^2/4$} & $c_2 \le c < 0$ \\ \cline{2-3}
\multicolumn{1}{|c|}{} & \multicolumn{1}{c|}{$a^2/4 < b \le a^2/3$} & $c_2 \le c \le c_1 < 0$  \\ \hline
\end{tabular}
\end{center}
\end{table}

\begin{table}[!h]
\begin{center}
\begin{tabular}{|ccc|}
\hline
\multicolumn{3}{|c|}{{\bf (II) Three Negative Roots}}  \\ \hline
\multicolumn{1}{|c|}{$\bm{a}$} & \multicolumn{1}{c|}{$\bm{b}$} & $\bm{c}$ \\ \hline
\multicolumn{1}{|c|}{\multirow{2}{*}{$a > 0$}} & \multicolumn{1}{c|}{$0 < b \le a^2/4$} & $0 < c \le c_1$ \\ \cline{2-3}
\multicolumn{1}{|c|}{} & \multicolumn{1}{c|}{$a^2/4 < b \le a^2/3$} & $0 < c_2 \le c \le c_1$  \\ \hline
\end{tabular}
\end{center}
\end{table}

\begin{table}[!h]
\begin{center}
\begin{tabular}{|ccc|}
\hline
\multicolumn{3}{|c|}{{\bf (III) Two Positive Roots and}}  \\
\multicolumn{3}{|c|}{{\bf One Negative Root}}  \\ \hline
\multicolumn{1}{|c|}{$\bm{a}$} & \multicolumn{1}{c|}{$\bm{b}$} & $\bm{c}$ \\ \hline
\multicolumn{1}{|c|}{\multirow{3}{*}{$a < 0$}} & \multicolumn{1}{c|}{$b < 0$} & $0 < c \le c_1$ \\ \cline{2-3}
\multicolumn{1}{|c|}{} & \multicolumn{1}{c|}{$b = 0$} & $0 < c \le -4a^3/27$ \\ \cline{2-3}
\multicolumn{1}{|c|}{} & \multicolumn{1}{c|}{$0 < b \le a^2/4$} & $0 < c \le c_1$ \\ \hline
\multicolumn{1}{|c|}{$a = 0$} & \multicolumn{1}{c|}{$b < 0$} & $-(2\sqrt{3}/9)\sqrt{-b^3} \le c < 0$ \\ \hline
\multicolumn{1}{|c|}{$a > 0$} & \multicolumn{1}{c|}{$b < -a^2/9$} & $0 < c \le c_1$ \\ \hline
\end{tabular}
\end{center}
\end{table}

\begin{table}[!h]
\begin{center}
\begin{tabular}{|ccc|}
\hline
\multicolumn{3}{|c|}{{\bf (IV) One Positive Root and }}  \\
\multicolumn{3}{|c|}{{\bf Two Negative Roots}}  \\ \hline
\multicolumn{1}{|c|}{$\bm{a}$} & \multicolumn{1}{c|}{$\bm{b}$} & $\bm{c}$ \\ \hline
\multicolumn{1}{|c|}{$a < 0$} & \multicolumn{1}{c|}{$b < -a^2/9$} & $c_2 \le c < 0$ \\ \hline
\multicolumn{1}{|c|}{$a = 0$} & \multicolumn{1}{c|}{$b < 0$} & $-(2\sqrt{3}/9)\sqrt{-b^3} \le c < 0$ \\ \hline
\multicolumn{1}{|c|}{\multirow{3}{*}{$a > 0$}} & \multicolumn{1}{c|}{$b < 0$} & $c_2 \le c < 0$ \\ \cline{2-3}
\multicolumn{1}{|c|}{} & \multicolumn{1}{c|}{$b = 0$} & $-4a^3/27 \le c < 0$ \\ \cline{2-3}
\multicolumn{1}{|c|}{} & \multicolumn{1}{c|}{$0 < b \le a^2/4$} & $c_2 \le c < 0$ \\ \hline
\end{tabular}
\end{center}
\end{table}


\begin{table}[!h]
\begin{center}
\begin{tabular}{|ccc|}
\hline
\multicolumn{3}{|c|}{{\bf (V) One Positive Root and }}  \\
\multicolumn{3}{|c|}{{\bf a Pair of Complex-Conjugate Roots}}  \\ \hline
\multicolumn{1}{|c|}{$\bm{a}$} & \multicolumn{1}{c|}{$\bm{b}$} & $\bm{c}$ \\ \hline
\multicolumn{1}{|c|}{\multirow{4}{*}{$a < 0$}} & \multicolumn{1}{c|}{$b < 0$} & $c < c_2 < 0$ \\ \cline{2-3}
\multicolumn{1}{|c|}{} & \multicolumn{1}{c|}{$b = 0$} & $c < 0$ \\ \cline{2-3}
\multicolumn{1}{|c|}{} & \multicolumn{1}{c|}{$0 < b \le a^2/3$} & $c < c_2 < 0$ \\ \cline{2-3}
\multicolumn{1}{|c|}{} & \multicolumn{1}{c|}{$b > a^2/3$} & $c < 0$ \\ \hline
\multicolumn{1}{|c|}{\multirow{2}{*}{$a = 0$}} & \multicolumn{1}{c|}{$b < 0$} &  $c < -(2\sqrt{3}/9)\sqrt{-b^3} < 0$ \\ \cline{2-3}
\multicolumn{1}{|c|}{} & \multicolumn{1}{c|}{$b \ge 0$} & $c < 0$ \\ \hline
\multicolumn{1}{|c|}{\multirow{4}{*}{$a > 0$}} & \multicolumn{1}{c|}{$b < 0$} & $c < c_2 < 0$ \\ \cline{2-3}
\multicolumn{1}{|c|}{} & \multicolumn{1}{c|}{$b = 0$} & $c < -4a^3/27 < 0$ \\ \cline{2-3}
\multicolumn{1}{|c|}{} & \multicolumn{1}{c|}{$0 < b \le a^2/4$} & $c < c_2 \le 0$ \\ \cline{2-3}
\multicolumn{1}{|c|}{} & \multicolumn{1}{c|}{$b > a^2/4$} & $c < 0$ \\ \hline
\end{tabular}
\end{center}
\end{table}


\newpage

\begin{table}[!h]
\begin{center}
\begin{tabular}{|ccc|}
\hline
\multicolumn{3}{|c|}{{\bf (VI) One Negative Root and }}  \\
\multicolumn{3}{|c|}{{\bf a Pair of Complex-Conjugate Roots}}  \\ \hline
\multicolumn{1}{|c|}{$\bm{a}$} & \multicolumn{1}{c|}{$\bm{b}$} & $\bm{c}$ \\ \hline
\multicolumn{1}{|c|}{\multirow{4}{*}{$a < 0$}} & \multicolumn{1}{c|}{$b < 0$} & $c > c_1 > 0$ \\ \cline{2-3}
\multicolumn{1}{|c|}{} & \multicolumn{1}{c|}{$b = 0$} & $c > -4a^3/27 > 0$ \\ \cline{2-3}
\multicolumn{1}{|c|}{} & \multicolumn{1}{c|}{$0 < b \le a^2/4$} & $c > c_1 \ge 0$ \\ \cline{2-3}
\multicolumn{1}{|c|}{} & \multicolumn{1}{c|}{$b > a^2/4$} & $c > 0$ \\ \hline
\multicolumn{1}{|c|}{\multirow{2}{*}{$a = 0$}} & \multicolumn{1}{c|}{$b < 0$} &  $c > (2\sqrt{3}/9)\sqrt{-b^3} > 0$ \\ \cline{2-3}
\multicolumn{1}{|c|}{} & \multicolumn{1}{c|}{$b \ge 0$} & $c > 0$ \\ \hline
\multicolumn{1}{|c|}{\multirow{4}{*}{$a > 0$}} & \multicolumn{1}{c|}{$b < 0$} & $c > c_1 > 0$ \\ \cline{2-3}
\multicolumn{1}{|c|}{} & \multicolumn{1}{c|}{$b = 0$} & $c < 0$ \\ \cline{2-3}
\multicolumn{1}{|c|}{} & \multicolumn{1}{c|}{$0 < b \le a^2/3$} & $c > c_1 > 0$ \\ \cline{2-3}
\multicolumn{1}{|c|}{} & \multicolumn{1}{c|}{$b > a^2/3$} & $c > 0$ \\ \hline
\end{tabular}
\end{center}
\end{table}

\section{Example}
\n
Consider the equation
\b
x^3 + 3 x^2 - \frac{1}{2} x - 4 = 0.
\e
One has: $a = 3$, $b = -1/2$ and $c = -4$.

\begin{itemize}
\item [{\bf (i)}] For the given $a = 3$ and $b =- 1/2$, one finds that $- a^2/9 < b < 0$ and, hence, $-c_1 < 0 < -c_0 < -ab  < -c_2$. The regime that applies is {\bf (2)} from Section 3.
\item[{\bf (ii)}] From (\ref{c0}) and (\ref{rho}) and from using $\rho_0 = - a/3$, one also finds:
\begin{itemize}
\item[{\bf (a)}] $- c_0 = 2.5000$,
\item[{\bf (b)}] $ - ab = 3.5000$,
\item[{\bf (c)}] $\rho_1 = 0.8708, \,\, \rho_2 = -2.8708$,
\item[{\bf (d)}] $\rho_0 = -1.0000$.
\end{itemize}
\item [{\bf (iii)}] The roots of quadratic equations (\ref{q1}), (\ref{q2}), and (\ref{q3}), together with $\xi_{1,2}$ from (\ref{xi}), are, respectively:
\begin{itemize}
\item[{\bf (a)}] $c_1 = 0.0203, \,\, c_2 = - 5.0203$. The given $c = -4$ is therefore between $c_2$ and $c_1$.
\item[{\bf (b)}] $\mu_1 = 0.0801, \,\, \mu_2 = - 2.0801$. Hence, the given cubic has either 1 or 3 real roots. In view of the fact that $c_2 < c < c_1$, the cubic has three real roots.
\item[{\bf (c)}] $\lambda_1 = 0.1583, \,\, \lambda_2 = - 3.1583$. The fact that $\lambda_{1,2}$ are real could be expected, as $c_1$ and $c_2$ have opposite signs. Nevertheless, one needs $\lambda_{1,2}$ as these are endpoints of some of the isolation intervals.
\item[{\bf (d)}] $\xi_1 = - 3.1602, \,\, \xi_2 = 1.1602.$
\end{itemize}
\item [{\bf (iv)}] The list $-c_2$, $-c_1$, $-c_0$, $-ab$, $0$ and the given $-c$ in increasing order is:
$-c_1 < 0 < -c_0 < -ab  < -c < -c_2$.
\item [{\bf (v)}] $a = 3$ and the equation is not depressed. Skip this point.
\item [{\bf (vi)}] Given that $a > 0$, one can see that Figure 5 is the relevant one.
\item [{\bf (vii)}] The neighbours of $-c = 4$ in the list from {\bf (iii)} are $-ab = 3.5000$ from below and $-c_2 = 5.0203$ from above. The intersection points of the straight lines $- b x - c_2 = (1/2)x + 5.0203$ and $- bx - ab = (1/2) x + 3.5000$ with $x^2(x + 3)$ allow the determination of the isolation intervals of the cubic $x^3 + 3 x^2 - (1/2) x - 4$. Namely, point (5) in the caption of Figure 5 applies and thus the isolation intervals of the roots are as follows:
    $\bm{-a \le x_3 \le  \mu_2}$, $\bm{\,\, \mu_2 \le x_2 \le -\sqrt{-b}}\,\,$ and $\bm{\,\,\sqrt{-b} \le x_1 \le \xi_2}$. \\
    Substituting the values of the endpoints of the isolation intervals, one finds: \\
    $-3 \le x_3 \le  - 2.0801$, $- 2.0801 \le x_2 \le -0.7071$, and $0.7071 \le x_1 \le 1.1602$.
\item [{\bf (viii)}] For the given value of $c$, one does not need to search for root bounds.
\item [{\bf (ix)}] The ``root harness" (\ref{harness}) is $\sqrt{3} \sqrt{a^2/3 - b} \le I \le 2 \sqrt{a^2/3 - b}$. For $a = 3$ and $b = -1/2$, one gets $3.2403 \le I \le 3.7417$. This can be narrowed down. Given that $ - ab < -c < -c_2$ (see Figure 5), the length of the interval which contains the three roots, that is  $x_1 + |x_3|$, decreases from $\sqrt{-b} + |a| = 3.7071$ (when $-c = - ab = 3/2$) to $\xi_2 + | \mu_2 | = 3.2403$ (when $-c = -c_2 = 5.0203$). Hence, $x_1 + | x_3 | $ cannot be smaller than $3.2403$ and cannot be bigger than $3.7071$. \\
    Solved with Maple, the roots of the equation $x^3 + 3 x^2 - (1/2) x - 4 = 0$ are: $x_3 = -2.6010, \,\, x_2 = -1.4556,$ and $x_1 = 1.0566$ --- each of which is exactly within its corresponding isolation interval, as determined above. \\
    Also, $x_1 + |x_3| = 3.6576$, which is within the narrowed ``root harness" $[3.2403,$ $3.7071 ]$. The cubic equation with $c = ab = -3/2$ has roots $x_3 = a = -3$, $x_2 = -\sqrt{-b} = -0.7071$, and $x_1 = \sqrt{-b} = 0.7071$ (Figure 5). These are spread over an interval of length $3.7071$ --- the upper bound of the ``root harness". The cubic equation with $c = c_2 = - 5.0203$ has the double root $x_{2,3} = \mu_2 = - 2.0801$ and the root $x_1 = \xi_2 = 1.1602$ (Figure 5). The distance between these roots is exactly $3.2403$ --- the lower bound of the ``root harness".
\end{itemize}

\section{Another Example --- Rayleigh Waves}

\n
The Rayleigh waves are elastic surface waves. These propagate on the surface of solids (decaying exponentially inside the solids). Some seismic waves are Rayleigh waves. A Rayleigh wave is a superposition of two waves propagated independently \cite{ll}. These are a combination of compression and dilation, resulting in elliptical motion. The equation of motion is \cite{ll}
\b
\label{eom}
\frac{\partial^2 \vec{u}}{\partial t^2} - c^2 \Delta \vec{u} = 0,
\e
where $\vec{u} = \vec{u}_l + \vec{u}_t$ is the displacement vector with $\vec{u}_l$ being the displacement in the direction of propagation (the co-called longitudinal wave) and $\vec{u}_t$ is the displacement in a plane, perpendicular to the direction of propagation (transverse wave) \cite{ll}. In (\ref{eom}), $c$ is either $c_l$ or $c_t$ --- the velocities of the two waves (often referred to as longitudinal and transverse speeds of sound) \cite{ll}. The relationship between these two speeds can be determined as follows. Consider the modulus of hydrostatic compression $K$ (also known as bulk modulus) and the modulus of rigidity $\mu$ (also known as shear modulus). Both $K$ and $\mu$ are always positive. These are connected \cite{ll} through the Lam\'e coefficient $\lambda$, that is: $\lambda = K - 2 \mu/3$. The Young modulus $E$ is given by \cite{ll}: $E = 9 K \mu/(3 K + \mu)$. The Poisson ratio $\sigma$ of transverse compression to longitudial dilation is given by \cite{ll}: $\sigma = (3K - 2 \mu)/(6 K + 2 \mu)$. In view of the positivity of the two moduli, $\sigma$ varies between $-1$ (when $K = 0$) and $1/2$ (when $\mu = 0$) \cite{ll}. One can express $K$ and $\mu$ from these two formul\ae \,\, \cite{ll}: $K = E/(3 - 6 \sigma)$ and $\mu = E/(2 + 2\sigma)$. The longitudinal and transverse speeds of sounds are \cite{ll}:
\b
c_l & = & \sqrt{\frac{3K + 4\mu}{3 \rho}} \,\,\, = \,\,\, \sqrt{\frac{\lambda + 2\mu}{\rho}} \,\,\, = \,\,\,  \sqrt{\frac{E (1 - \sigma)}{\rho (1 + \sigma)(1 - 2\sigma)}}, \\
c_t & = & \sqrt{\frac{\mu}{\rho}} \,\,\, = \,\,\,  \sqrt{\frac{E}{2 \rho (1 + \sigma)}},
\e
where $\rho$ is the density. Hence, $c_l > \sqrt{4/3} c_t$. \\
The stress-tensor is given by \cite{ll}:
\b
\sigma_{ik} = \frac{E}{1 + \sigma} \left( u_{ik} + \frac{\sigma}{1 - 2 \sigma} u_{ll} \delta_{ik} \right)
\e
(sum over the repeated index $l$). In the above, $i, k, l$ are $x, y, z$ and $u_{ik}$ are the components of the strain tensor (it is symmetric) \cite{ll}:
\b
u_{ik} = \frac{1}{2} \left( \frac{\partial u_i}{\partial x_k} + \frac{\partial u_k}{\partial x_i} + \frac{\partial u_l}{\partial x_i} \frac{\partial u_l}{\partial x_k} \right)
\e
(sum over the repeated index $l$). Here $u_i$ are the components of the displacement vector $\vec{u} = \vec{r}\phantom{'}' - \vec{r}$, with $\vec{r}$ and $\vec{r}\phantom{'}'$ --- the radius-vectors of a point before and after the deformation, respectively. \\
Suppose the surface of the elastic medium is the $xy$-plane and that the medium is in $z < 0$. Hence, $\sigma_{xz} = \sigma_{yz} = \sigma_{zz} = 0$ and the conditions
$u_{xz} = 0, \,\, u_{yz} = 0, \,\, \sigma(u_{xx} + u_{yy}) + (1 - \sigma) u_{zz} = 0$ ensue \cite{ll}. \\
The transverse part of the wave, $\vec{u}_t$, is divergence-free: div $\vec{u}_t = 0$, while the longitudinal part of the wave, $\vec{u}_l$, is rotation-free: curl $\vec{u}_l = 0$  \cite{ll}. These conditions lead to:
\b
a (k^2 - \kappa_t^2) + 2 b k \kappa_l & = & 0, \\
2 a k \kappa_t + b (k^2 + \kappa_t^2) & = & 0,
\e
where $a$ and $b$ are some constants, $k$ is the wave number, $\kappa_t = \sqrt{k^2 - \omega^2/c_t^2}$ and $\kappa_l = \sqrt{k^2 - \omega^2/c_l^2}$ are the rapidities of the transverse and longitudinal damping, respectively, and $\omega$ is the angular frequency  \cite{ll}. The condition for compatibility of these two is given by \cite{ll}:
\b
\label{cond}
\left( 2 k^2 - \frac{\omega^2}{c_t^2} \right)^4 = 16 k^4 \left( k^2 - \frac{\omega^2}{c_t^2} \right) \left( k^2 - \frac{\omega^2}{c_l^2} \right).
\e
Introducing $\xi = \omega/(c_t k) > 0$, $x = \xi^2$, and $q = c_t^2/c_l^2$, from (\ref{cond}), one gets the cubic equation
\b
\label{landau}
x^3 - 8x^2 + 8 (3 - 2q) x - 16(1 -q) = 0.
\e
Clearly, the allowed values of the parameter $q$ satisfy $0 \le q < 3/4$. \\
Additionally, the right-hand side of (\ref{cond}) must be non-negative, as is the left-hand side. This leads to $(1 + \xi)(1- \xi)(1 - \sqrt{q} \xi)(1 + \sqrt{q} \xi) \ge 0$. Hence, either $0 < \xi \le 1$ or $\xi \ge 1/\sqrt{q}$ --- only such roots are physical. Depending on $q$, the latter is greater than or equal to $\sqrt{4/3} \approx 1.1547$.  \\
To align with the notation used, introduce $a = -8$, $b = 8(3 - 2q)$ (note that $b$ is linear in $q$ and is always positive for $0 \le q < 3/4$), and $c = -16(1 - q)$ (also linear in $q$ and always negative for $0 \le q < 3/4$). Equation (\ref{landau}) then takes the form $x^3 + ax^2 + bx + c = 0$. \\
The roots (\ref{c12}) of the first auxiliary quadratic equation (\ref{q1}) are:
\b
c_{1,2} = c_0 \,\, \pm \,\, \frac{2}{27} \sqrt{(a^2 - 3b)^3} = -\frac{704}{27} + \frac{128}{3} q \,\, \pm \,\, \frac{16 \sqrt{8}}{27} \sqrt{(6 q - 1)^3}.
\e
Hence, if $0 \le q < 1/6$, the cubic (\ref{landau}) will have one real root only. For $1/6 \le q < 3/4$, the cubic will have three real roots if the coefficient $c = 16 q - 16$ is between $c_2$ and $c_1$. Indeed, $c > c_2$ always, while $c_1 - c$ changes sign (from negative to positive) when $q$ varies from $0$ to $3/4$. Hence, there is a value of $q$, above which the equation (\ref{landau}) has 3 real roots. \\
Next, one needs to determine the end-points of the isolation intervals. The roots (\ref{mu12}) of the second auxiliary quadratic equation (\ref{q2}) are
\b
\mu_{1,2} = \frac{8}{3} \pm \frac{2\sqrt{2}}{3} \sqrt{6q - 1}.
\e
Both $\mu_{1,2}$ are positive for $1/6 \le q < 4/3$. \\
The corresponding simple roots (\ref{xi}) are:
\b
\xi_{1,2} = \frac{8}{3} \pm \frac{4\sqrt{2}}{3} \sqrt{6q - 1}.
\e
The roots of the cubic with $c = c_0 = -2a^3/27 + ab/3 = -704/27 + (128/3) q$ are:
\b
\rho_0 & = & - \frac{a}{3} \,\,\, = \,\,\, 8/3, \\
\rho_{1,2} & = & \frac{8}{3} \pm \frac{2\sqrt{6}}{3} \sqrt{6q - 1}.
\e
Finally, the non-zero roots of the separatrix equation $x^3 + a x^2 + bx = 0$ are:
\b
\lambda_{1,2} = 4 \pm 2 \sqrt{2} \sqrt{2q-1}.
\e
Note that $-c/b = (2q - 2)/(2q-3)$. This is $2/3$ when $q = 0$ and decreases monotonically towards $1/3$ when $q \to 3/4$. \\
Note that $0 < b \le 2a^2/9 = 128/9$ for $11/18 \le  q < 3/4$. In this range of $q$, one has $c_2 < 0$, $c_1 > 0$, $c < 0$, and $-c_1 < 0 < -c < -c_2$. Thus Figure 10 (4) applies. Equation (\ref{landau}) has three positive roots in the case of $11/18 <  q < 3/4$ and their isolation intervals are:
\begin{itemize}
\item[{\bf (ii)}] $(2q - 2)/(2q-3) < x_3 <  8/3 - (2\sqrt{2}/3) \sqrt{6q - 1},$
\item[{\bf (ii)}] $8/3 - (2\sqrt{2}/3) \sqrt{6q - 1} < x_2 <  4 - 2 \sqrt{2} \sqrt{2q - 1},$
\item[{\bf (iii)}] $4 - 2 \sqrt{2} \sqrt{2q - 1} < x_1 < 8/3 - (2\sqrt{6}/3) \sqrt{6q - 1}.$
\end{itemize}
Next, for $1/2 \le q < 11/18$, one has $2a^2/9 < b \le a^2/4$. In this range of $q$, one again has $c_2 < 0$, $c_1 > 0$, $c_0 < 0$, $c < 0$, and $-c_1 < 0 -c_0 < -c < -c_2$. Hence, Figure 12(4) is the relevant one and equation (\ref{landau}) has three positive roots in the case of $1/2 <  q < 11/18$ with isolation intervals given by:
\begin{itemize}
\item[{\bf (i)}] max$\{-8/3 - (2\sqrt{6}/3) \sqrt{6q-1}, \,\,\, (2q - 2)/(2q-3) \} < x_3 < 8/3 - (2\sqrt{2}/3) \sqrt{6q - 1},$
\item[{\bf (ii)}] $8/3 - (2\sqrt{2}/3) \sqrt{6q - 1} < x_2 < 8/3,$
\item[{\bf (iii)}] $8/3 + (2\sqrt{6}/3) \sqrt{6q-1} < x_1 < 8/3 - (4\sqrt{2}/3) \sqrt{6q-1}.$
\end{itemize}
Next, for $1/6 \le q < 1/2$, one has $a^2/4 < b \le a^2/3$. In this range of $q$, one has $c_2 < 0$, $c_1 < 0$ (note that $c_1 = 0$ when $q = 1/2$), $c_0 < 0$, and $c < 0$.
While $-c < -c_2$ for $1/6 \le q < 1/2$, there is a value of $q$ in this range, say $\tilde{q}$, below which $-c < -c_1$. This $\tilde{q}$ is the only real root of the cubic equation $c_1 = c$ and will not be determined. Also, at $q = 17/45 > \tilde{q}$, one has $c = c_0$ and for values of $q$ below $17/45$, one has $-c < -c_0$. Hence, for $17/45 < q < 1/2$, one has $0 < -c_1 < -c_0 < -c < -c_2$ with Figure 14(4) applying. The isolation intervals of the three positive roots for $17/45 < q < 1/2$ are the same as those for the case of $1/2 <  q < 11/8$ above. \\
For values of $q$ satisfying $\tilde{q} < q < 17/45$, one has $0 < -c_1 < -c < -c_0 < -c_2$. Figure 14(3) applies. There are three positive roots with the following isolation intervals:
\begin{itemize}
\item[{\bf (i)}] max$\{8/3 + (4\sqrt{2}/3) \sqrt{6q-1}, \,\,\, (2q - 2)/(2q-3)  \} < x_3 <  8/3 - (2\sqrt{6}/3) \sqrt{6q-1},$
\item[{\bf (ii)}] $ 8/3 < x_2 < -8/3 + (2\sqrt{2}/3) \sqrt{6q-1},$
\item[{\bf (iii)}] $8/3 + (2\sqrt{2}/3) \sqrt{6q-1} < x_1 < 8/3 + (2\sqrt{6}/3) \sqrt{6q-1}.$
\end{itemize}
For values of $q$ such that $1/6 < q < \tilde{q}$, one has $0 < -c < -c_1 < -c_0 < -c_2$. Figure 14(2) applies. The cubic equation (\ref{landau}) has a single real root (positive) with isolation interval given by $(2q - 2)/(2q-3) < x_1 <  8/3 + (4\sqrt{2}/3) \sqrt{6q-1}$. \\
Finally, for $0 \le q < 1/6$, one has $b > a^2/3$. For this range of $q$, $c_{1,2}$ are not real. One further has $c_0 < 0$, $c < 0$, $ab < 0$, and $0 < -c < -c_0 < -ab$. Hence, when $0 \le q < 1/6$, Figure 16(2) applies. Equation (\ref{landau}) has only one real root (positive). The isolation interval of this root is $(2q - 2)/(2q-3) < x_1 < 8/3$.


\begin{thebibliography}{99}

\bibitem{card} G. Cardano, {\it Ars Magna or the Rules of Algebra, Translated and Edited by T. Richard Witmer}, Dover (2000).

\bibitem{g} R. Garver, {\it A Reading List in the Elementary Theory of Equations}, The American Mathematical Monthly {\bf 40(2)}, 77--84 (1933) doi: 10.2307/2300939

\bibitem{rs} Q.I. Rahman and G. Schmeisser, {\it Analytic Theory of Polynomials}, Oxford University Press (2002).

\bibitem{p} V.V. Prasolov, {\it Polynomials}, Springer (2010).

\bibitem{m} M. Marden, {\it Geometry of Polynomials, Math. Surveys no. 3}, American Mathematical Society, Providence, RI (1966).

\bibitem{d} L.E. Dickson, {\it First Course in the Theory of Equations}, Braunworth (1922).

\bibitem{k} A. Kurosh, {\it Higher Algebra}, Mir Publishers (1980).

\bibitem{c} R. Descartes, {\it La G\'eom\'etrie} (1637).

\bibitem{ga} C.F. Gauss, {\it Werke, Dritter Band}, G\"ottingen (1876), p.67.

\bibitem{newton} I. Newton, {\it Universal Arithmetick or, A Treatise of Arithmetical Composition and Resolution}, W. Johnston London (1769).

\bibitem{st} J.C.F. Sturm, {\it M\'emoire sur la R\'esolution des \'Equations Num\'eriques}, Bulletin G\'en\'eral et Universel des Annonces et des Nouvelles Scientifiques (Bulletin des Sciences de F\'erussac), {\bf 11}, 419--425 (1829).

\bibitem{bu} F.D. Budan, {\it Nouvelle M\'ethode pour la R\'esolution des \'Equations Num\'eriques}, Paris: Courcier (1807).

\bibitem{cheng} S.S. Cheng and Y.-Z. Lin, {\it Dual Sets of Envelopes and Characteristic Regions of Quasi-Polynomials}, World Scientific (2009).

\bibitem{a}  D.S. Arnon, {\it Geometric Reasoning with Logic and Algebra}, Artificial Intelligence {\bf 37}, 37--60 (1988).

\bibitem{3} E.M. Prodanov, {\it Classification of the Roots of the Quartic Equation and their Pythagorean Tunes}, International Journal of Applied and Computational Mathematics (Springer) {\bf 7}, 218 (2021), doi: 10.1007/s40819-021-01152-w, arXiv: 2008.07529.

\bibitem{28} E.M. Prodanov, {\it A Method for Locating the Real Roots of the Symbolic Quintic Equation Using Quadratic Equations}, Advanced Theory and Simulations (Wiley) (2022), 2200011, doi: 10.1002/adts.202200011, arXiv:2106.02977.

\bibitem{1} E.M. Prodanov, {\it On the Determination of the Number of Positive and Negative Polynomial Zeros and Their Isolation}, Open Mathematics (de Gruyter) {\bf 18}, 1387--1412 (2020), doi: 10.1515/math-2020-0079, arXiv: 1901.05960.

\bibitem{4} E.M. Prodanov, {\it The Siebeck--Marden--Northshield Theorem and the Real Roots of the Symbolic Cubic Equation}, Resultate der Mathematik (Birkha\"user) {\bf 77}, 126 (2022), doi: 10.1007/s00025-022-01667-8, arXiv:2107.01847.

\bibitem{2} E.M. Prodanov, {\it New Bounds on the Real Polynomial Roots}, Comptes Rendus de l'Acad\'emie Bulgare des Sciences {\bf 75(2)} 178--186 (2022), doi: 10.7546/CRABS.2022.02.02 , arXiv:2008.11039.

\bibitem{ll} L.D. Landau and E.M. Lifshitz, {\it Theory of Elasticity (Course of Theoretical Physics, Volume 7)}, Pergamon Press (1970).

\end{thebibliography}
\end{document}